\documentclass[12pt]{amsart}

\usepackage{amsmath}
\usepackage{amssymb}
\usepackage{geometry}
\usepackage{enumitem}
\usepackage{extpfeil}
\usepackage[all,2cell]{xy} \UseAllTwocells

\usepackage{hyperref}
\usepackage{xcolor}
\definecolor{darkgreen}{rgb}{0,0.5,0}
\hypersetup{colorlinks,urlcolor=blue,citecolor=darkgreen,linkcolor=blue,linktocpage=false}

\usepackage[abbrev]{amsrefs}
\usepackage[capitalise,noabbrev]{cleveref}

\numberwithin{subsection}{section}
\numberwithin{equation}{section}

\newtheorem{theorem}[equation]{Theorem}
\newtheorem{proposition}[equation]{Proposition}
\newtheorem{corollary}[equation]{Corollary}
\newtheorem{lemma}[equation]{Lemma}

\newtheorem{innercustomthm}{Theorem}
\newenvironment{customthm}[1]
  {\renewcommand\theinnercustomthm{#1}\innercustomthm}
  {\endinnercustomthm}

\theoremstyle{remark}
\newtheorem{remark}[equation]{Remark}
\newtheorem{example}[equation]{Example}
\newtheorem{warning}[equation]{Warning}

\theoremstyle{definition}
\newtheorem{assumption}[equation]{Assumption}
\newtheorem{definition}[equation]{Definition}
\newtheorem{notation}[equation]{Notation}

\newextarrow{\xrightrightarrow}{{20}{20}{20}{20}}{\bigRelbar\bigRelbar{\bigtwoarrowsleft\rightarrow\rightarrow}}

\newcommand{\F}{\mathbb{F}}
\newcommand{\G}{\mathbb{Z}/p^2}

\newcommand{\Z}{\mathbb{Z}}

\newcommand{\cat}[1]{\mathcal{#1}} %
\newcommand{\steen}{\mathcal{A}}
\newcommand{\TT}{\mathcal{T}}

\newcommand{\al}{\alpha}
\newcommand{\be}{\beta}

\newcommand{\de}{\delta}
\newcommand{\del}{\partial}
\newcommand{\ep}{\epsilon}
\newcommand{\Ga}{\Gamma}

\newcommand{\io}{\iota}

\newcommand{\Si}{\Sigma}
\newcommand{\si}{\sigma}
\newcommand{\te}{\theta}

\newcommand{\bu}{\bullet}
\newcommand{\pb}{\ar@{}[dr]|{\mbox{\huge $\lrcorner$}}}
\newcommand{\po}{\ar@{}[dr]|{\mbox{\huge $\ulcorner$}}}

\newcommand{\ral}{\xrightarrow} %
\newcommand{\Ra}{\Rightarrow}
\newcommand{\surj}{\twoheadrightarrow}
\newcommand{\tild}{\widetilde}
\newcommand{\tworal}{\xrightrightarrow} %

\newcommand{\dfn}{:=}
\newcommand{\inv}{\boxminus}
\newcommand{\lan}{\left\langle}

\newcommand{\op}{\oplus}
\newcommand{\ot}{\otimes}
\newcommand{\ran}{\right\rangle}
\newcommand{\sm}{\wedge}
\newcommand{\sq}{\square \,}

\newcommand{\x}{\times}

\newcommand{\Ch}{\mathbf{Ch}}
\newcommand{\ct}[1]{\mathbf{#1}} %
\newcommand{\denorm}{\Theta} %
\newcommand{\Gpd}{\mathbf{Gpd}}
\newcommand{\Mod}{\mathbf{Mod}}
\newcommand{\Topp}{\mathbf{Top}}

\DeclareMathOperator{\Aut}{Aut}
\DeclareMathOperator{\coker}{coker}

\DeclareMathOperator{\Hom}{Hom}
\DeclareMathOperator{\im}{im}
\DeclareMathOperator{\Mon}{Mon}
\DeclareMathOperator{\Ob}{Ob}
\DeclareMathOperator{\tr}{Tr}

\newcommand{\ab}{\mathrm{ab}}
\newcommand{\id}{\mathrm{id}}
\newcommand{\proj}{\mathrm{proj}}
\newcommand{\sh}{\mathrm{sh}}
\newcommand{\strict}{\mathrm{strict}}

\newcommand{\Def}[1]{\textbf{\boldmath{#1}}} 

\newdir{->}{!/5pt/\dir{-}!/2.5pt/\dir{-}*!/-5pt/\dir2{>}}

\begin{document}

\title{The DG-category of secondary cohomology operations} 

\author{Hans-Joachim Baues}
\email{baues@mpim-bonn.mpg.de}
\address{Max-Planck-Institut f\"ur Mathematik\\
Vivatsgasse 7\\
53111 Bonn\\
Germany}

\author{Martin Frankland}
\email{Martin.Frankland@uregina.ca}
\address{University of Regina\\
3737 Wascana Parkway\\
Regina, Saskatchewan, S4S 0A2\\
Canada}

\date{\today}


\begin{abstract}
We study track categories (i.e., groupoid-enriched categories) endowed with additive structure similar to that of a $1$-truncated DG-category, except that composition is not assumed right linear. We show that if such a track category is right linear up to suitably coherent correction tracks, then it is weakly equivalent to a $1$-truncated DG-category. This generalizes work of the first author on the strictification of secondary cohomology operations. As an application, we show that the secondary integral Steenrod algebra is strictifiable.
\end{abstract}

\maketitle

\section{Introduction}

Cohomology operations are important tools in algebraic topology. The Steenrod algebra (of primary stable mod $p$ cohomology operations) was determined as a Hopf algebra in celebrated work of Milnor \cite{Milnor58}. The structure of secondary cohomology operations was determined as a ``secondary Hopf algebra'' in \cite{Baues06}, and via different methods in \cite{Nassau12}. Unlike for primary operations, composition of secondary operations is not bilinear, but bilinear up to homotopy. Part of the work in \cite{Baues06} was to \emph{strictify} the structure of secondary operations, i.e., replace it with a weakly equivalent differential bigraded algebra, in which composition is bilinear. The purpose of this paper is to revisit this strictification step, simplify it, and generalize it. 

Here is the motivating example in more detail. For a fixed prime number $p$, mod $p$ cohomology operations correspond to maps between finite products of Eilenberg--MacLane spaces $K(\F_p, n)$, the representing objects. Stable operations correspond to maps between finite products of Eilenberg--MacLane spectra $\Si^n H\F_p$. Primary operations are encoded by homotopy classes of such maps. More precisely, the Steenrod algebra $\steen$ is given by homotopy classes of maps
\[
\steen^n = [H\F_p, \Si^n H\F_p].
\]
For higher order cohomology operations, one needs more than homotopy classes. One way to encode higher order operations is the topologically enriched category $\cat{EM}$ consisting of finite products of Eilenberg--MacLane spectra
\[
A = \Si^{n_1} H\F_p \x \ldots \x \Si^{n_k} H\F_p
\]
and mapping \emph{spaces} between them. Composition in the homotopy category $\pi_0 \cat{EM}$ is bilinear, but composition in $\cat{EM}$ is not bilinear. It is strictly left linear, i.e., satisfies $(a+b)x = ax + bx$, and right linear up to coherent homotopy $a(x+y) \sim ax + ay$. The higher coherence for right linearity is studied in \cite{BauesFrankland17}.

For secondary operations, it suffices to take the fundamental groupoid of each mapping space in $\cat{EM}$. This yields a track category $\Pi_1 \cat{EM}$, i.e., a category enriched in groupoids. In fact, $\Pi_1 \cat{EM}$ has some additional additive structure. Since each mapping space in $\cat{EM}$ is an abelian group object, the same is true of $\Pi_1 \cat{EM}$. Now, an abelian group object in groupoids corresponds to a $1$-truncated chain complex $C_1 \to C_0$. Moreover, composition in $\Pi_1 \cat{EM}$ is left linear (strictly) and right linear up to track. Hence, the track category $\Pi_1 \cat{EM}$ looks like a $1$-truncated DG-category (i.e., a category enriched in $1$-truncated chain complexes), except that composition is not right linear. One of the structural results from \cite{Baues06} is the following.

\begin{customthm}{A}\label{thm:StrictificationSteenrod}
The track category $\Pi_1 \cat{EM}$ is weakly equivalent to a $1$-truncated DG-category over $\Z/p^2$.
\end{customthm}

The proof relied on correction tracks for right linearity $a(x+y) \Ra ax + ay$. These linearity tracks can be chosen to satisfy certain coherence conditions, which we call the \emph{linearity track equations}. The main result of this paper is the following; see \cref{thm:Strictification}.

\begin{customthm}{B}[Strictification theorem] 
Let $\cat{T}$ be a left linear track category which admits linearity tracks satisfying the linearity track equations. Then $\cat{T}$ is weakly equivalent to a $1$-truncated DG-category.

If moreover every morphism in $\cat{T}$ is $p$-torsion (i.e., satisfies $px = 0$), then $\cat{T}$ is weakly equivalent to a $1$-truncated DG-category over $\Z/p^2$.
\end{customthm}

The contribution of this paper is threefold.
\begin{itemize}
\item We streamline the construction of the strictification, which is not about secondary cohomology operations, but rather about coherence in track categories. This part is mostly expository, to make the relevant literature more transparent. Moreover, the current presentation can be adapted to tertiary cohomology operations.
\item One new result is the observation that the construction works over $\Z$, i.e., without \linebreak $p$-torsion assumption (Proposition~\ref{pr:PseudoFunctorZ}). As an application, we show that the secondary integral Steenrod algebra is strictifiable (Corollary~\ref{cor:IntegralEM}).
\item We provide an alternate proof of the strictification theorem using a $2$-categorical observation due to Lack. This bypasses the cocycle computation in Baues--Wirsching cohomology, the argument used in \cite{Baues06}.
\end{itemize}

\subsection*{Organization}

In \cref{sec:LinearityTracks}, we describe the notion of a track category $\cat{F}$ having some additive structure that makes composition left linear (strictly) and right linear up to coherent homotopy (Definition~\ref{def:LinTrackEq}). In \cref{sec:LeftLinTrack}, we show that such a track category in which composition is also right linear (strictly) is the same as a $1$-truncated DG-category. 
Next, the proof of the strictification theorem consists of three steps.

\begin{itemize}
\item The construction of a certain pseudo-functor $s \colon \cat{B}_0 \to \cat{F}$. This is done in \cref{sec:PseudoFunctor}.
\item Upgrading this construction to a certain pseudo-functor $s \colon \cat{B} \to \cat{F}$, where $\cat{B}$ is a $1$-truncated DG-category. This is done is \cref{sec:PseudoFunctorAction}.
\item Some general categorical facts about pseudo-functors ensuring that we obtain the desired weak equivalence. This is done in \cref{sec:Strictification}.
\end{itemize}

\cref{sec:ChoiceGraph} makes the general construction more explicit in the case of secondary cohomology operations. \cref{sec:Toda} explains how a strictification of $\cat{T}$ can be used to compute Toda brackets in $\cat{T}$.

\subsection*{Related work}

There are other strictification problems in track categories with additive structure. The \emph{strengthening theorem} \cite{BauesJP03}*{Theorem~6.2.2} says that under certain assumptions, a track category with weak products is weakly equivalent to a track category with strict products. If the track category has weak products and weak coproducts, then one cannot in general strictify both the products and coproducts simultaneously. Gaudens showed that one can strictify the products and make the weak coproducts somewhat more strict \cite{Gaudens10}.

Using Baues--Wirsching cohomology of small categories along with calculations in \linebreak Hochschild, Shukla, and MacLane cohomology, the first author and Pirashvili recovered the strictification theorem for the secondary Steenrod algebra (\cref{thm:StrictificationSteenrod}) and generalized it \cite{BauesP04}*{Theorem~8.1.1} \cite{BauesP06}, cf.\ \cite{BauesJP08}*{\S 3}. The current paper makes no use of cohomology theories for categories. It is not obvious whether one could prove the strictification theorem for the secondary integral Steenrod algebra (Corollary~\ref{cor:IntegralEM}) using a similar cohomological argument.

There is also literature on the strictification of pseudo-algebras for certain $2$-monads on certain $2$-categories \cite{Power89} \cite{Lack02cod} \cite{Shulman12}. It would be interesting to see if left linear track categories equipped with linearity tracks form the pseudo-algebras of some appropriate $2$-monad whose strict algebras are the bilinear track categories.

\subsection*{Acknowledgements}

We thank David Blanc, Teimuraz Pirashvili, and Emily Riehl for helpful conversations. We also thank the referee for their careful reading and useful suggestions. The second author thanks the Max-Planck-Institut f\"ur Mathematik Bonn for its generous hospitality. The second author was partially funded by a grant of the Deutsche Forschungsgemeinschaft SPP~1786: Homotopy Theory and Algebraic Geometry, as well as the Natural Sciences and Engineering Research Council of Canada (NSERC), Discovery Grant RGPIN-2019-06082. Cette recherche a \'et\'e financ\'ee par le Conseil de recherches en sciences naturelles et en g\'enie du Canada (CRSNG), subvention D\'ecouverte RGPIN-2019-06082.

\section{Preliminaries and notation}

\begin{notation}
A \Def{groupoid} is a category in which every morphism is invertible. 
In this paper, we will only consider \emph{small} groupoids. 
Denote the data of a (small) groupoid by $G = \left( G_0, G_1, \de_0, \de_1, \id^{\sq}, \sq, (-)^{\inv} \right)$, where:
\begin{itemize}
 \item $G_0 = \Ob(G)$ is the set of objects of $G$.
 \item $G_1 = \Hom(G)$ is the set of morphisms of $G$, also called \Def{tracks} in $G$. The set of morphisms from $x$ to $y$ is denoted $G(x,y)$. We consider a groupoid $G$ as a graded set, with
\[
\deg(x) = \begin{cases}
0 &\text{if } x \in G_{0} \\
1 &\text{if } x \in G_{1} \\
\end{cases}
\]
and we write $x \in G$ in each case.
 \item $\de_0 \colon G_1 \to G_0$ is the source map.
 \item $\de_1 \colon G_1 \to G_0$ is the target map.
 \item $\id^{\sq} \colon G_0 \to G_1$ sends each object $x$ to its corresponding identity morphism $\id^{\sq}_{x}$.
 \item $\sq \colon G_1 \x_{G_0} G_1 \to G_1$ is composition in $G$.
 \item $f^{\inv} \colon y \to x$ is the inverse of the morphism $f \colon x \to y$.
\end{itemize}
Groupoids form a category $\Gpd$, where morphisms are functors between groupoids.

Denote the fundamental groupoid of a topological space $X$ by $\Pi_{1}(X)$.
\end{notation}

\begin{notation}
A groupoid $G$ is \Def{pointed} if it is equipped with a morphism of groupoids $\ast \to G$ from the terminal groupoid $\ast$ (with one object and one morphism). Let $\Gpd_*$ denote the category of pointed groupoids.
\end{notation}

\begin{definition}
A \Def{track category} is a category enriched in $(\Gpd,\x)$, the category of groupoids with its Cartesian product as monoidal structure.

A track category $\TT$ is \Def{pointed} if every mapping groupoid $\TT(A,B)$ is pointed, with basepoint denoted $0 = 0_{A,B} \in \TT(A,B)_0$, and for any objects $A,B,C$ of $\TT$, the composition map
\[
\mu \colon \TT(B,C) \x \TT(A,B) \to \TT(A,C)
\]
satisfies the following two conditions:
\begin{itemize}
\item Objects: $\mu(0,y) = 0$ and $\mu(x,0) = 0$ for all $x \in \TT(B,C)_0$ and $y \in \TT(A,B)_0$.
\item Morphisms: $\mu(\id^{\sq}_0,b) = \id^{\sq}_0$ and $\mu(a,\id^{\sq}_0) = \id^{\sq}_0$ for all $a \in \TT(B,C)_1$ and $b \in \TT(A,B)_1$.
\end{itemize}
By abuse of notation, we will sometimes write $0 \in \TT(A,B)_1$ for $\id^{\sq}_0$.
\end{definition}

\begin{remark}
A pointed track category is the same as a category enriched in $(\Gpd_*,\sm)$. 
The smash product of pointed groupoids $G \sm G' = (G \x G')/(G \vee G')$ makes $(\Gpd_*, \sm)$ into a symmetric monoidal category \cite{ElmendorfM09}*{Lemma~4.20}.
\end{remark}

\begin{definition}
The \Def{homotopy category} of a track category $\cat{T}$ is the category $\pi_0 \cat{T}$ with the same objects as $\cat{T}$ and whose hom-sets are obtained by taking components of each mapping groupoid:
\[
\left( \pi_0 \cat{T} \right)(A,B) = \pi_0 \left( \cat{T}(A,B) \right).
\]
The \Def{underlying category} of $\cat{T}$ is the category $\cat{T}_0$ obtained by forgetting the tracks, i.e., with hom-sets $\cat{T}_0(A,B) = \cat{T}(A,B)_0$.
\end{definition}

We write $x \in \cat{T}$ if $x \in \cat{T}(A,B)$ for some objects $A$ and $B$. For $x,y \in \cat{C}$, we write $xy = \mu(x,y)$ when $x$ and $y$ are composable, i.e., when the target of $y$ is the source of $x$, and $\deg(x) = \deg(y)$ holds. From now on, whenever an expression such as $xy$ appears, it is understood that $x$ and $y$ must be composable.

\begin{notation}
For $\deg(x) = \deg(y) = 0$ and $\deg(\al) = 1$, denote:
\begin{equation} \label{eq:TensorFromPointwise}
\begin{cases}
x \ot y \dfn xy \\
x \ot \al \dfn \id^{\sq}_x \al \, \text{, also written } x \al \\
\al \ot x \dfn \al \id^{\sq}_x \, \text{, also written } \al x. \\
\end{cases}
\end{equation}
We call $x \ot y$ the \Def{$\ot$-composition} of $x$ and $y$, which is defined whenever $\deg(x) + \deg(y) \leq 1$ holds. The $\ot$-composition is associative, unital, and satisfies $\deg(x \ot y) = \deg(x) + \deg(y)$. Moreover, it determines the pointwise composition. Indeed, for $\deg(\al) = \deg(\be) = 1$, the following factorizations hold in $\cat{T}$:
\begin{align}
\al \be &= \left( \al \ot \de_1 \be \right) \sq \left( \de_0 \al \ot \be \right) \nonumber \\  
&= \left( \de_1 \al \ot \be \right) \sq \left( \al \ot \de_0 \be \right). \label{eq:PointwiseFromTensor}
\end{align}
For our purposes, it will be more convenient to work with the $\ot$-composition instead of the pointwise composition.
\end{notation}

If $\cat{T}$ is a pointed track category and $\al$ and $\be$ are tracks to zero, i.e., satisfying $\de_1 \al = 0$ and $\de_1 \be = 0$, then Equation~\eqref{eq:PointwiseFromTensor} specializes to
\begin{equation}\label{eq:FactorizToZero}
(\de_0 \al) \ot \be = \al \ot (\de_0 \be).
\end{equation}

\section{Linearity tracks}\label{sec:LinearityTracks}

The purpose of this paper is to study distributivity in track categories, i.e., the compatibility between multiplicative and additive structure. In this section, we describe the additive structure of interest, where composition is left linear (strictly) and right linear up to coherent homotopy.

\begin{definition}\label{def:LeftLinTrack}
A \Def{locally linear track category} $\cat{T}$ is a pointed track category such that each mapping groupoid $\cat{T}(A,B)$ is an abelian group object in $\Gpd$ (based at $0_{A,B}$). 
The track category $\cat{T}$ is \Def{left linear} if moreover composition in $\cat{T}$ is left linear, i.e., satisfies $(a+a')x = ax + a'x$. \Def{Right linear} and \Def{bilinear} are defined analogously.

A \Def{morphism} of locally linear track categories is a track functor (i.e., $\Gpd$-enriched functor) $F \colon \cat{S} \to \cat{T}$ such that for all objects $A,B$ of $\cat{S}$, the induced map of groupoids $F \colon \cat{S}(A,B) \to \cat{T}(FA,FB)$ is a map of abelian group objects in groupoids, i.e., preserves addition (strictly). 
\end{definition}

\begin{lemma}\label{lem:LeftLinearity}
For $\cat{T}$ a locally linear track category, the following are equivalent.
\begin{enumerate}
\item Composition in $\cat{T}$ is left linear.
\item The abelian group object structure of $\cat{T}(A,B)$ is natural in $A$.
\end{enumerate}
If moreover $\cat{T}$ has finite products (in the $\Gpd$-enriched sense, namely strict products), 
then these conditions are further equivalent %
to the following.
\begin{enumerate}
\setcounter{enumi}{2}
\item For every object $B$, the abelian group object structure of $\cat{T}(A,B)$ is induced by a (strict) 
abelian group object structure on $B$, i.e., by ``pointwise addition in the target''.
\end{enumerate}
\end{lemma}

\begin{definition}
Let $\cat{T}$ be a left linear track category, and let $x,y \colon X \to A$ and $a \colon A \to B$ be maps in $\cat{T}$. A track $\Ga_a^{x,y} \in \cat{T}(X,B)_1$ of the form 
\[
\Ga_a^{x,y} \colon a(x+y) \Ra ax + ay
\]
is called a \Def{linearity track}.
\end{definition}

\begin{definition} \label{def:LinTrackEq}
The following are called the \Def{linearity track equations}.
\begin{enumerate}
\item \label{item:TrackPrecompo} \emph{Precomposition}:
$\Ga_a^{xz,yz} = \Ga_a^{x,y} z$. In other words, the following diagram of tracks commutes:
\[
\xymatrix{
a (xz + yz) \ar@{=}[d] \ar@{=>}[r]^-{\Ga_a^{xz,yz}} & axz + ayz \ar@{=}[d] \\
a(x+y)z \ar@{=>}[r]_-{\Ga_a^{x,y} z} & (ax + ay)z. \\
}
\]
\item \emph{Postcomposition}:
$\Ga_{ba}^{x,y} = \Ga_b^{ax,ay} \sq b \Ga_a^{x,y}$. As a diagram:
\[
\xymatrix{
ba (x + y) \ar@{=>}[d]_{b \Ga_a^{x,y}} \ar@{=>}[r]^-{\Ga_{ba}^{x,y}} & bax + bay \\
b(ax+ay). \ar@{=>}[ur]_-{\Ga_b^{ax,ay}} & \\
}
\]
In particular, setting $a=b=1$ yields the unital equation $\Ga_1^{x,y} = \id^{\sq}_{x+y}$.
\item \emph{Symmetry}:
$\Ga_a^{x,y} = \Ga_a^{y,x}$. 
\item \emph{Left linearity}:
$\Ga_{a+a'}^{x,y} = \Ga_a^{x,y} + \Ga_{a'}^{x,y}$. 
\item \label{item:TrackAssoc} \emph{Associativity}: $\left( \Ga_a^{x,y} + a z \right) \sq \Ga_{a}^{x+y, z} = \left( a x + \Ga_{a}^{y,z} \right) \sq \Ga_{a}^{x,y+z}$. As a diagram:
\[
\xymatrix{
a (x+y+z) \ar@{=>}[d]_{\Ga_a^{x,y+z}} \ar@{=>}[r]^-{\Ga_a^{x+y,z}} & a (x+y) + az \ar@{=>}[d]^{\Ga_a^{x,y} + az} \\
ax + a(y+z) \ar@{=>}[r]_-{ax + \Ga_a^{y,z}} & ax + ay + az. \\
}
\]
In particular, setting $y=z=0$ yields $\Ga_a^{x,0} = \id^{\sq}_{ax}$ and likewise $\Ga_a^{0,y} = \id^{\sq}_{ay}$.
\item \label{item:NaturalXY} \emph{Naturality in $x$ and $y$}: Given tracks $G \colon x \Ra x'$ and $H \colon y \Ra y'$, the equation
\[
\left( aG + aH \right) \sq \Ga_a^{x,y} = \Ga_{a}^{x',y'} \sq a (G + H) 
\]
holds in $\cat{T}_1$. As a diagram:
\[
\xymatrix{
a(x+y) \ar@{=>}[d]_{a(G+H)} \ar@{=>}[r]^-{\Ga_a^{x,y}} & ax + ay \ar@{=>}[d]^{aG + aH} \\
a(x'+y') \ar@{=>}[r]_-{\Ga_a^{x',y'}} & ax' + ay'. \\
}
\]
\item \label{item:NaturalA} \emph{Naturality in $a$}: Given a track $\al \colon a \Ra a'$, the equation
\[
\left( \al x + \al y \right) \sq \Ga_{a}^{x,y} = \Ga_{a'}^{x,y} \sq \al (x+y)
\]
holds in $\cat{T}_1$. As a diagram:
\[
\xymatrix{
a(x+y) \ar@{=>}[d]_{\al(x+y)} \ar@{=>}[r]^-{\Ga_a^{x,y}} & ax + ay \ar@{=>}[d]^{\al x + \al y} \\
a'(x+y) \ar@{=>}[r]_-{\Ga_{a'}^{x,y}} & a'x + a'y. \\
}
\]
\end{enumerate}
\end{definition}

Let us recall how linearity tracks arise \cite{Baues06}*{\S 4.2}.

\begin{proposition} \label{pr:LinTrackEq}
Let $\cat{T}$ be a left linear track category with finite (strict) products. Assume that for every object $A$ of $\cat{T}$, the two inclusion maps $i_1 = (1_A, 0) \colon A \to A \x A$ and $i_2 = (0, 1_A) \colon A \to A \x A$ exhibit $A \x A$ as a weak coproduct, i.e., the restriction
\[
\xymatrix{
\cat{T}(A \x A, B) \ar[r]^-{(i_1^*,i_2^*)}_-{\sim} & \cat{T}(A,B) \x \cat{T}(A,B) \\
}
\]
is an equivalence of groupoids for every object $B$ of $\cat{T}$. Then $\cat{T}$ admits canonical linearity tracks, which moreover satisfy the linearity track equations.
\end{proposition}

\begin{proof}
For every map $a \colon A \to B$, let $\Ga_a \in \cat{T}(A \x A, B)_1$ be the unique track satisfying the equations
\[
\begin{cases}
i_1^* \Ga_a = \id^{\sq}_a \\
i_2^* \Ga_a = \id^{\sq}_a. \\
\end{cases}
\]
For every $x,y \colon X \to A$, define the composite $\Ga_a^{x,y} := \Ga_a \ot (x,y) \in \cat{T}(X,B)_1$, which is a track $\Ga_a^{x,y} \colon a(x+y) \Ra ax + ay$ as illustrated in the diagram
\[
\xymatrix{
X \ar[r]^-{(x,y)} \ar[dr]_-{x+y} & A \x A \ar[d]_{+_A} \ar[r]^-{a \x a} & B \x B \ar[d]^{+_B} \\
& A \ar[r]^-{a} & B. \ultwocell<\omit>{\Ga_a} \\
}
\]
These tracks $\Ga_a^{x,y}$ satisfy the linearity track equations \cite{Baues06}*{Theorem~4.2.5}.
\end{proof}

Note that for such a track category $\cat{T}$, the homotopy category $\pi_0 \cat{T}$ is additive.

\begin{example}
If $\cat{C}$ is a topologically enriched category satisfying the topological analogue of the assumptions of Proposition~\ref{pr:LinTrackEq}, then the proposition applies to the underlying track category $\cat{T} = \Pi_1 \cat{C}$. This happens in the example of higher order cohomology operations, as described in \cite{BauesFrankland17}.
\end{example}

\subsection{Iterated linearity tracks}

For the remainder of the section, let $\cat{T}$ be left linear track category equipped with system of linearity tracks $\Ga_a^{x,y} \colon a(x+y) \Ra ax + ay$ satisfying the linearity track equations.

\begin{definition}
Given an integer $n \geq 2$ and maps $x_1, \ldots, x_n \colon X \to A$ and $a \colon A \to B$, define the track $\Ga_a^{x_1, \ldots, x_n} \colon a(x_1 + \ldots + x_n) \Ra ax_1 + \ldots + ax_n$ inductively by
\[
\Ga_a^{x_1, \ldots, x_n} \dfn \left( \Ga_{a}^{x_1, \ldots , x_{n-1}} + ax_n \right) \sq \Ga_a^{x_1 + \ldots + x_{n-1}, x_n}
\]
as illustrated in the diagram
\[
\xymatrix @C=3.5pc {
a \left( x_1 + \ldots + x_{n-1} + x_n \right) \ar@{=>}[d]_{\Ga_a^{x_1 + \ldots + x_{n-1}, x_n}} \ar@{=>}[r]^-{\Ga_a^{x_1, \ldots, x_n}} & ax_1 + \ldots + ax_n \\
a \left( x_1 + \ldots + x_{n-1} \right) + ax_n. \ar@{=>}[ur]_-{\Ga_a^{x_1 + \ldots + x_{n-1}} + ax_n} & \\
}
\]
For $n=1$, take by convention the identity track $\Ga_a^{x_1} = \id^{\sq} \colon ax_1 \Ra ax_1$.
\end{definition}

\begin{proposition}\label{pr:BreakSum}
\begin{enumerate}
\item The $(n-1)!$ ways of distributing the product $a(x_1 + \ldots + x_n)$ via a complete bracketing, counted by ordering the $n-1$ instances of the symbol~$+$, all yield the same track $\Ga_a^{x_1, \ldots, x_n} \colon a(x_1 + \ldots + x_n) \Ra ax_1 + \ldots + ax_n$. 
\item Writing the sum into $k$ blocks
\begin{align*}
x_1 + \ldots + x_n &= \overbrace{\left( x_1 + \ldots + x_{n_1} \right)}^{n_1} + \overbrace{\left( x_{n_1 + 1} + \ldots + x_{n_1 + n_2} \right)}^{n_2} + \ldots + \overbrace{\left( x_{n_1 + \ldots + n_{k-1} + 1} + \ldots + x_{n} \right)}^{n_k} \\
&=: S_1 + \ldots + S_k
\end{align*}
yields the factorization
\[
\Ga_a^{x_1, \ldots, x_n} = \left( \sum_{i=1}^k \Ga_a^{x_{n_1 + \ldots + n_{i-1} + 1}, \ldots, x_{n_1 + \ldots + n_{i}}} \right) \sq \Ga_a^{S_1, \ldots, S_k}
\]
as illustrated in the diagram
\[
\xymatrix @C-3pc {
**[l] a \left( x_1 + \ldots + x_n \right) = a (S_1 + \ldots + S_k) \ar@{=>}[dr]_-{\Ga_a^{x_1, \ldots, x_n}} \ar@{=>}[r]^-{\Ga_a^{S_1, \ldots, S_k}} & a S_1 + \ldots + a S_k \ar@{=>}[d]^-{\Ga_a^{x_1, \ldots, x_{n_1}} + \ldots + \Ga_a^{x_{n_1 + \ldots + n_{k-1} + 1}, \ldots, x_n}} \\
& ax_1 + \ldots + ax_n. \\
}
\]
\end{enumerate}
\end{proposition}

\begin{proof}
The case $n=3$ holds by assumption, as Equation~\ref{pr:LinTrackEq}~\eqref{item:TrackAssoc}. The case $n=4$ says that the diagram
\[
\xymatrix @C-0.8pc {
& a(w+x) + a(y+z) \ar@{=>}[dd]_-(0.7){\Ga_a^{w,x} + a(y+z)} \ar@{=>}[rr]^-{a(w+x) + \Ga_a^{y,z}} & & a(w+x) + ay + az \ar@{=>}[dd]_-{\Ga_a^{w,x} + ay + az} \\
a(w+x+y+z) \ar@{=>}[dd]_-{\Ga_a^{w,x+y+z}} \ar@{=>}[ur]^-{\Ga_{a}^{w+x,y+z}} \ar@{=>}[rr]^-(0.32){\Ga_a^{w+x+y,z}} & & a(w+x+y) + az \ar@{=>}[dd]_-(0.7){\Ga_a^{w,x+y} + az} \ar@{=>}[ur]^-{\Ga_{a}^{w+x,y} +az} & \\
& aw + ax + a(y+z) \ar@{=>}[rr]^-(0.3){aw + ax + \Ga_a^{y,z}} & & aw + ax + ay + az \\
aw + a(x+y+z) \ar@{=>}[ur]^-{aw + \Ga_{a}^{x,y+z}} \ar@{=>}[rr]^-{aw + \Ga_a^{x+y,z}} & & aw + a(x+y) + az \ar@{=>}[ur]^-{aw + \Ga_{a}^{x,y} + az} & \\
}
\]
commutes. The front face commutes by induction, and is equal to $\Ga_a^{w,x+y,z}$; likewise for the top and left faces. The right face commutes by induction, and is equal to $\Ga_a^{w,x,y} + az$; likewise for the bottom face. The back face commutes and is equal to $\Ga_a^{w,x} + \Ga_a^{y,z}$, by the interchange law in the additive groupoid $\cat{T}(X,B)$. The general case $n \geq 4$ is proved similarly by induction.

The second statement is a straightforward generalization of the factorization of $\Ga_a^{w,x,y,z}$ :
\[
\xymatrix @C=3.5pc {
a(w+x+y+z) \ar@{=>}[dr]_-{\Ga_a^{w,x,y,z}} \ar@{=>}[r]^-{\Ga_a^{w+x, y+z}} & a(w+x) + a(y+z) \ar@{=>}[d]^-{\Ga_a^{w,x} + \Ga_a^{y,z}} \\
& aw + ax + ay + az \\
}
\]
using the back face of the cube.
\end{proof}

\begin{proposition} \label{pr:IteratedLinTrackEq}
The iterated linearity tracks $\Ga_{a}^{x_1, \ldots, x_n}$ satisfy the following equations, which are analogous to the linearity track equations in Proposition~\ref{pr:LinTrackEq}.
\begin{enumerate}
\item \emph{Precomposition}:
$\Ga_a^{x_1 z, \ldots, x_n z} = \Ga_a^{x_1, \ldots, x_n} z$. 
\item \emph{Postcomposition}:
$\Ga_{ba}^{x_1, \ldots, x_n} = \Ga_b^{ax_1, \ldots, ax_n} \sq b \Ga_a^{x_1, \ldots, x_n}$. 
\item \label{item:IteratedSymm} \emph{Symmetry}:
$\Ga_a^{x_1, \ldots, x_n} = \Ga_a^{x_{\si(1)}, \ldots, x_{\si(n)}}$ for any permutation $\si \in \Si_n$. 
\item \emph{Left linearity}:
$\Ga_{a+a'}^{x_1, \ldots, x_n} = \Ga_a^{x_1, \ldots, x_n} + \Ga_{a'}^{x_1, \ldots, x_n}$. 
\item \emph{Naturality in $x_i$}: Given tracks $G_i \colon x_i \Ra x_i'$ for $1 \leq i \leq n$, the equation
\[
\left( aG_1 + \ldots + aG_n \right) \sq \Ga_a^{x_1, \ldots, x_n} = \Ga_{a}^{x'_1, \ldots, x'_n} \sq a (G_1 + \ldots + G_n) 
\]
holds. 
\item \emph{Naturality in $a$}: Given a track $\al \colon a \Ra a'$, the equation
\[
\left( \al x_1 + \ldots + \al x_n \right) \sq \Ga_{a}^{x_1, \ldots x_n} = \Ga_{a'}^{x_1, \ldots, x_n} \sq \al (x_1 + \ldots + x_n)
\]
holds. 
\end{enumerate}
\end{proposition}

\begin{proof}
This follows inductively from the case $n=2$.
\end{proof}

\subsection{Multiplying by an integer}

\begin{notation}
For a map $a \colon A \to B$ in $\cat{T}_0$ and $n \geq 1$, denote the track in $\cat{T}(A,B)_1$
\[
\Ga(n)_a \dfn \Ga_a^{1_A, \ldots, 1_A} \colon a(\overbrace{1_A + \ldots + 1_A}^{n \text{ terms}}) \Ra a 1_A + \ldots + a 1_A = n \cdot a.
\]
In particular, if $p$ annihilates every map in $\cat{T}$ and $p \vert n$, then $\Ga(n)$ is a track of the form $\Ga(n) \colon 0 \Ra 0$. 
\end{notation}

\begin{remark}\label{rem:pTorsion}
The $p$-torsion condition $px =0$ is meant for morphisms $x \in \cat{T}$ of degree $0$, but together with left linearity, this implies that tracks are also $p$-torsion. Indeed, let $\al \colon x \Ra y$ a be track between morphisms $x,y \colon A \to B$. Then we have
\[
p \al = \al + \ldots + \al = (1_B + \ldots + 1_B) \al = 0.
\]
\end{remark}

\begin{lemma} \label{lem:GammaP2}
For a map $a \colon A \to B$ in $\cat{T}_0$ and $m,n \geq 1$, the following equality holds:
\[
\Ga(m \cdot n)_a = (m \cdot 1_B) \Ga(n)_a \sq \Ga(m)_a (n \cdot 1_A)
\]
in $\cat{T}(A,B)_1$. In other words, the following diagram of tracks commutes:
\[
\xymatrix @C=4pc {
a (mn \cdot 1_A) \ar@{=}[d] \ar@{=>}[rr]^-{\Ga(mn)_a} & & mn \cdot a \ar@{=}[d] \\
a (m \cdot 1_A) (n \cdot 1_A) \ar@{=>}[r]_-{\Ga(m)_a (n \cdot 1_A)} & (m \cdot a) (n \cdot 1_A) = (m \cdot 1_B) a (n \cdot 1_A) \ar@{=>}[r]_-{(m \cdot 1_B) \Ga(n)_a} & (m \cdot 1_B) (n \cdot a). \\
}
\]
In particular, if $p$ annihilates every map in $\cat{T}_0$ and $p^2 \vert n$, then we have $\Ga(n) = \id^{\sq}_0 \colon 0 \Ra 0$.
\end{lemma}

\begin{proof}
Break the sum $mn \cdot 1_A$ into $m$ blocks of $n$ terms each
\[
mn \cdot 1_A = \overbrace{(n \cdot 1_A) + \ldots + (n \cdot 1_A).}^{m \text{ blocks}}
\]
Using Proposition~\ref{pr:BreakSum}, we obtain
\begin{flalign*}
&&\Ga(mn)_a &= \Ga_a^{\overbrace{1_A, \ldots, 1_A}^{mn}} && \\
&& &= \left( \sum_1^m \Ga_a^{\overbrace{1_A, \ldots, 1_A}^{n}} \right) \sq \Ga_a^{\overbrace{n \cdot 1_A, \ldots, n \cdot 1_A}^{m}} && \\
&& &= m \cdot \Ga(n)_a \sq \Ga_a^{\overbrace{n \cdot 1_A, \ldots, n \cdot 1_A}^{m}} && \\
&& &= m \cdot \Ga(n)_a \sq \Ga_a^{\overbrace{1_A, \ldots, 1_A}^{m}} (n \cdot 1_A) &&\text{by Equation } \ref{pr:LinTrackEq} \eqref{item:TrackPrecompo} \\
&& &= m \cdot \Ga(n)_a \sq \Ga(m)_a (n \cdot 1_A). &&
\end{flalign*}
\end{proof}

Next, we deal with negatives.

\begin{notation}
For any map $a \colon A \to B$ in $\cat{T}_0$, define the track $\Ga(-1)_a \colon a(- 1_A) \Ra -a$ by the commutative diagram of tracks
\[
\xymatrix{
a(1 + (-1)) \ar@{=}[d] \ar@{=>}[rr]^-{\Ga_a^{1,-1}} & & a(1) + a(-1) \ar@{=>}[d]^{a + \Ga(-1)_a} \\
a(0) \ar@{=}[r] & 0 \ar@{=}[r] & a + (-a). \\
}
\]
Explicitly, it is given by $\Ga(-1)_a = (\Ga_a^{1,-1} - a)^{\inv} = - \Ga_a^{1,-1} + a(-1)$.
\end{notation}
The analogously defined track $a(-x) \Ra -ax$ for an arbitrary map $x \colon X \to A$ is equal to $\Ga(-1)_a x$, by the precomposition equation. 

\begin{lemma} \label{lem:MultiplyNegative}
For any map $a \in \cat{T}_0$ and integer $m > 0$, the following diagram of tracks commutes:
\[
\xymatrix{
a(-m) \ar@{=>}[d]_{\Ga_a^{-1,\ldots,-1}} \ar@{=>}[r]^-{\Ga(-1)_a (m)} & - a(m) \ar@{=>}[d]^{- \Ga(m)_a} \\
a(-1) + \ldots + a(-1) \ar@{=>}[r]^-{m \Ga(-1)_a} & -m a. \\
}
\]
Denote the resulting track by $\Ga(-m)_a \colon a(-m) \Ra -m a$.
\end{lemma}

\begin{lemma} \label{lem:SumNegatives}
For any map $a \in \cat{T}_0$ and integers $m,n \in \Z$, the following equality of tracks holds:
\[
\Ga(m+n)_a = (\Ga(m)_a + \Ga(n)_a) \sq \Ga_a^{m,n}.
\]
In diagrams:
\[
\xymatrix{
a(m+n) \ar@{=>}[dr]_{\Ga(m+n)_a}  \ar@{=>}[r]^-{\Ga_a^{m,n}} & a(m) + a(n) \ar@{=>}[d]^{\Ga(m)_a + \Ga(n)_a} \\
& ma + na = (m+n)a. \\
}
\]
\end{lemma}

\begin{proof}
The case $m,n \geq 0$ follows from Proposition~\ref{pr:BreakSum}. The general case $m,n \in \Z$ follows from Lemma~\ref{lem:MultiplyNegative}.
\end{proof}

\section{Left linear track categories and DG-categories} \label{sec:LeftLinTrack}

In this section, we consider pointed track categories endowed with a certain additive structure. The motivational example is when $\cat{C}$ is a category enriched in $(\Topp_*, \sm)$, and each mapping space $\cat{C}(A,B)$ has the structure of a topological abelian group. Note that $\cat{C}$ is \emph{not} enriched in topological abelian groups, as we do not assume that composition is bilinear. However, we will assume that composition is left linear, i.e., satisfies $(x+x')y = xy + x'y$, as is the case when addition of maps $x,x' \in \cat{C}(B,C)$ is defined pointwise in the target. We are interested in the $\Gpd_*$-category $\Pi_{1} \cat{C}$ of such a category $\cat{C}$.

\subsection{Abelian group objects in groupoids}

The following equivalence can be found in \cite{Baues06}*{Proposition~2.2.6}, \cite{Bourn07}*{Theorem~1.2, Remark~1}, or \cite{Bourn90}*{\S 2}. Here we fix some choices.

\begin{proposition} \label{pr:AbGpObjGpd}
The category $\Gpd_{\ab}$ of abelian group objects in the category $\Gpd$ of small groupoids is equivalent to the category of $1$-truncated chain complexes of abelian groups (in other words, chain complexes concentrated in degrees $0$ and $1$). The equivalence sends an abelian group object $G$ in $\Gpd$ to its \emph{Moore chain complex}
\[
M(G) \dfn \left( \ker \de_1 \ral{\del = \de_0} G_0 \right).
\]
An inverse equivalence assigns to a $1$-truncated chain complex of abelian groups $F_1 \ral{\del} F_0$ the groupoid denoted
\[
\denorm(F) \dfn \left( F_1 \op F_0 \tworal[\de_1]{\de_0} F_0 \right)
\]
defined as follows. For $(x_1,x_0) \in F_1 \op F_0$, the source and target maps are given by
\[
\begin{cases}
\de_0 (x_1, x_0) = \del x_1 + x_0 \\
\de_1 (x_1, x_0) = x_0 \\
\end{cases}
\]
so that $(x_1,x_0) \colon \del x_1 + x_0 \Ra x_0$ is a track in the groupoid $G(\del)$. The composition of tracks is given by
\begin{equation} \label{eq:CompoTracksAdd}
(x_1, x_0) \sq (y_1, y_0) = (x_1 + y_1, x_0)
\end{equation}
when the composability condition $\de_1 y = y_0 = \de_0 x = \del x_1 + x_0$ is satisfied. 

Likewise, the category of $\F_p$-vector space objects in $\Gpd$ is equivalent to the category of $1$-truncated chain complexes of $\F_p$-vector spaces.

The homotopy groups of the groupoid $\denorm(F)$ are given by the homology of the corresponding chain complex:
\[
\pi_i \denorm(F) \cong H_i (F) = \begin{cases}
\coker \del &\text{if } i=0 \\
\ker \del &\text{if } i=1 \\
0 &\text{otherwise.} \\
\end{cases}
\] 
\end{proposition}

Via the equivalence of Proposition~\ref{pr:AbGpObjGpd}, a left linear track category (as in Definition~\ref{def:LeftLinTrack}) can be viewed as the data $\cat{F} = \left( \cat{F}_1 \ral{\del} \cat{F}_0, +, \ot \right)$, where we replace each mapping groupoid $\cat{T}(A,B) = \left( \cat{T}(A,B)_1 \tworal[\de_1]{\de_0} \cat{T}(A,B)_0 \right)$ by the corresponding $1$-truncated chain complex of abelian groups
\[
M \cat{T}(A,B) = \cat{F}(A,B) = \left( \cat{F}(A,B)_1 \ral{\del} \cat{F}(A,B)_0 \right).
\]

\subsection{Truncated chain complexes}

In this section, a \emph{chain complex} will mean a non-negatively graded chain complex unless otherwise noted, i.e., a chain complex $C$ satisfying $C_i = 0$ for $i < 0$.  We work in the category $\Mod_R$ of $R$-modules, for some commutative ring $R$. The tensor product $C \ot D$ of chain complexes of $R$-modules will mean the tensor product $C \ot_R D$ over $R$ unless otherwise noted.

Let us recall some basics about truncation of chain complexes.

\begin{definition}
Let $n \geq 0$ be an integer.
\begin{enumerate}
\item A chain complex $C$ is called \Def{$n$-truncated} if it is trivial above degree $n$, that is, satisfying $C_i = 0$ for $i > n$. Denote by $\Ch_{\leq n}$ the full subcategory of $n$-truncated chain complexes and by $\io \colon \Ch_{\leq n} \to \Ch$ its inclusion into the category of all chain complexes.
\item The \Def{$n$-truncation} of a chain complex $C$ is the $n$-truncated chain complex
\[
(\tr_n C)_i = \begin{cases}
C_i &\text{if } i < n \\
\coker \left( C_{n+1} \ral{d} C_n \right) = C_n / \im d &\text{if } i = n \\
0 &\text{if } i > n \\
\end{cases}
\]
with differential inherited from that of $C$. This construction defines a functor \linebreak $\tr_n \colon \Ch \to \Ch_{\leq n}$. 
\end{enumerate}
\end{definition}

Recall that a (non-negatively graded) \Def{differential graded category}, or \Def{DG-category} for short, is a category enriched in chain complexes $(\Ch, \ot, R)$. 

\begin{definition}
A DG-category $\cat{F}$ is called \Def{$n$-truncated} if every hom-complex $\cat{F}(X,Y)$ is $n$-truncated. Note that this is the same as a category enriched in $\Ch_{\leq n}$, where the tensor product in $\Ch_{\leq n}$ is given by $M \ot_{n} N \dfn \tr_n(M \ot N)$.
\end{definition}

The $n$-truncation $\tr_n \colon \Ch \to \Ch_{\leq n}$ is also known as the \emph{good} $n$-truncation, because it induces the $n$-truncation on homology groups:
\[
H_i (\tr_n C) = \begin{cases}
H_i C &\text{if } i \leq n \\
0 &\text{if } i > n. \\
\end{cases}
\]
Moreover, $\tr_n$ is left adjoint to the inclusion, and the adjunction $\tr_n \dashv \io$ is monoidal. 

\begin{example}
A $0$-truncated DG-category over the ring $R = \Z$ is precisely a preadditive category. More generally, it is an $R$-linear category, i.e., a category enriched in $(\Mod_R, \ot_R)$.
\end{example}

\begin{example}\label{ex:1Alg}
Let us spell out explicitly the structure found in a $1$-truncated DG-category $\cat{F} = \left( \cat{F}_1 \ral{\del} \cat{F}_0, +, \ot \right)$.

\begin{enumerate}
\item A class of objects for $\cat{F}_0$.
\item For all objects $A$ and $B$ of $\cat{F}_0$, a $1$-truncated chain complex of $R$-modules
\[
\cat{F}(A,B) = \left( \cat{F}(A,B)_1 \ral{\del} \cat{F}(A,B)_0 \right).
\]
The zero elements are denoted $0 = 0_{A,B}\in \cat{F}(A,B)_0$.
\item For all object $A$, a distinguished unit element $1_A \in \cat{F}(A,A)_0$.
\item For $x,y \in \cat{F}$ composable and satisfying $\deg(x) + \deg(y) \leq 1$, the $\ot$-composition $x \ot y \in \cat{F}$ is defined and satisfies $\deg(x \ot y) = \deg(x) + \deg(y)$.
\end{enumerate}

The following equations are required to hold.

\begin{enumerate}
\item (\textit{Associativity}) $\ot$ is associative: $(x \ot y) \ot z = x \ot (y \ot z)$.
\item (\textit{Units}) The unit elements satisfy $x \ot 1 = x = 1 \ot x$ for all $x \in \cat{F}$.
\item (\textit{Bilinearity}) $\ot$ is bilinear.
\item (\textit{Leibniz rule}) The $\ot$-composition is a chain map, which yields the following equations. For $x, y, a, b \in \cat{F}$ with $\deg(x) = \deg(y) = 0$ and $\deg(a) = \deg(b) = 1$, we have:
\[
\begin{cases}
(\del a) \ot b = a \ot (\del b) \in \cat{F}_1 \\
\del (x \ot b) = x \ot (\del b) \in \cat{F}_0 \\
\del (a \ot y) = (\del a) \ot y \in \cat{F}_0. \\
\end{cases}
\]
\end{enumerate}
\end{example}

\begin{proposition} \label{pr:1Truncated}
A left linear track category which is right linear can be identified with a $1$-truncated DG-category (up to a strict track equivalence which is the identity on objects).
\end{proposition}

\begin{proof}
Like the Dold--Kan correspondence, the equivalence $M \colon \Gpd_{\ab} \cong \Ch_{\leq 1} \colon \denorm$ from Proposition~\ref{pr:AbGpObjGpd} is not a monoidal equivalence \cite{SchwedeS03equ}*{\S 2.3}. Both functors $M$ and $\denorm$ are lax monoidal, so that they induce change-of-enrichment functors \cite{Borceux94v2}*{Proposition~6.4.3}. 
The counit $\ep \colon M \denorm \ral{\cong} 1$ is monoidal, and the unit $\eta \colon 1 \ral{\cong} \denorm M$ is pseudo-monoidal. Nonetheless, applying the unit $\eta$ to each hom-groupoid of a bilinear track category $\cat{T}$ yields a pseudo-functor (see Definition~\ref{def:PseudoFunctor}) $\eta \colon \cat{T} \to \denorm M \cat{T}$ which turns out to be a (strict) track functor.
\end{proof}

\section{Construction of the pseudo-functor}\label{sec:PseudoFunctor}

In this section, starting from a left linear track category $\cat{F}$ which is right linear up to suitably coherent linearity tracks, we use the linearity tracks to construct a pseudo-functor $s \colon \cat{B}_0 \to \cat{F}$ that will later provide the desired strictification of $\cat{F}$.

\begin{definition}\label{def:PseudoFunctor}
Let $\cat{T}$ be a track category and $\cat{B}_0$ a category. A \Def{pseudo-functor} $(s,\Ga) \colon \cat{B}_0 \to \cat{T}$ consists of the following data.
\begin{enumerate}
\item A function assigning to each object $A$ of $\cat{B}_0$ an object $sA$ of $\cat{T}$.
\item For all objects $A$ and $B$ of $\cat{B}_0$, a function
\[
s \colon \cat{B}_0(A,B) \to \cat{T}(sA,sB)_0.
\]
\item For every (composable) $x,y \in \cat{B}_0$, a track
\[
\Ga(x,y) \colon (sx)(sy) \Ra s(xy).
\]
\end{enumerate} 
The following equations are required to hold.
\begin{enumerate}
\item (\textit{Associativity}) 
For every $x,y,z \in \cat{B}_0$, we have the equality
\begin{equation} \label{eq:CoherentIso}
\Ga(xy,z) \sq \left( \Ga(x,y) (sz) \right) = \Ga(x,yz) \sq \left( (sx) \Ga(y,z) \right)
\end{equation}
of tracks $(sx)(sy)(sz) \Ra s(xyz)$, as illustrated in the diagram
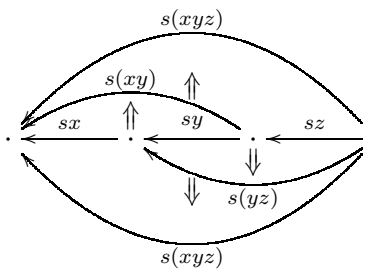
\[
\xymatrix @C+1pc {
\cdot & \cdot \ar[l]_{sx} & \cdot \lllowertwocell_{s(xy)}{} \ar[l]_{sy} & \cdot \ar[l]_{sz} \lluppertwocell^{s(yz)}{^} \xlowertwocell[lll]{}<-16>_{s(xyz)} \xuppertwocell[lll]{}<16>^{s(xyz)}{^} \\
}
\]
where pasting the four tracks yields the identity track $\id^{\sq}_{s(xyz)} \in \cat{T}_1$.
\item (\textit{Units}) For every object $A$ of $\cat{B}_0$ the equality $s(1_A) = 1_{sA}$ holds (strictly). For every $x \in \cat{B}_0(A,B)$, we have equalities $\Ga(1_B,x) = \id^{\sq}_{sx}$ and $\Ga(x,1_A) = \id^{\sq}_{sx}$ as tracks $sx \Ra sx$ in $\cat{T}_1$.
\end{enumerate}
\end{definition}

\begin{remark}
A pseudo-functor satisfying the strict unital condition above is sometimes called \emph{reduced}. This condition can be weakened to having tracks $1_{sA} \Ra s(1_A)$ that satisfy certain coherence conditions; cf.\ \cite{BauesM07}*{Appendix} and \cite{Borceux94v1}*{\S 7.5}. Our example of interest will satisfy the strict unital condition.
\end{remark}

As before, we fix a prime number $p$ and denote by $\F_p$ the field of $p$ elements. Consider the ring $\Z/p^2$ with the canonical quotient map $\G \surj \F_p$. Let $\cat{F}$ be a left linear track category in which every morphism is $p$-torsion, equipped with a system of linearity tracks $\Ga_a^{x,y} \colon a(x+y) \Ra ax + ay$ satisfying the linearity track equations. In this section, we construct a pseudo-functor
\[
(s,\Ga) \colon \cat{B}_0 \to \cat{F}
\]
which will induce a strictification of $\cat{F}$, as discussed in Section~\ref{sec:Strictification}. First, let us fix some notation and terminology.

\begin{notation}
A (directed) \Def{graph} $E = \left( E_0, E_1, \de_0, \de_1 \right)$ consists of sets $E_0$ and $E_1$, called the vertices and edges respectively, and two functions $\de_0, \de_1 \colon E_1 \to E_0$, called the source and target maps. A small category $\cat{C}$ has in particular an underlying graph $U \cat{C}$, and the forgetful functor $U \colon \ct{Cat} \to \ct{Graph}$ has a left adjoint
\[
\Mon \colon \ct{Graph} \to \ct{Cat}.
\]
We call $\Mon(E)$ the \Def{free category} generated by the graph $E$; cf.\ \cite{DwyerK80sim}*{\S 2}.
\end{notation}

Explicitly, the objects of $\Mon(E)$ are the vertices $E_0$ of $E$, and morphisms in $\Mon(E)$ are composable words in $E_1$. If $E_0 = \{ \ast \}$ consists of a single vertex, then $\Mon(E)$ is the free monoid on the set of edges $E_1$.

\begin{notation}
Given a commutative ring $R$ and a category $\cat{C}$, let $R \cat{C}$ denote the category with the same objects as $\cat{C}$, with morphisms modules in $R \cat{C}$ given by free $R$-modules
\[
\left( R \cat{C} \right)(A,B) := R \left( \cat{C}(A,B) \right) 
\]
and composition given by the $R$-bilinear extension of composition in $\cat{C}$, as illustrated in the diagram
\[
\xymatrix{
R \cat{C}(B,C) \ot_R R \cat{C}(A,B) \ar[d]_{\cong} \ar[r] & R \cat{C}(A,C) \\
R \left( \cat{C}(B,C) \x \cat{C}(A,B) \right). \ar[ur]_{R \mu} & \\ 
}
\]
\end{notation}

Now, choose a graph $E$ together with a graph morphism $h_E \colon E \to U \pi_0 \cat{F}$. Since the homotopy category $\pi_0 \cat{F}$ is $\F_p$-linear, this defines by adjunction an $\F_p$-linear functor \linebreak $h_E'' \colon \F_p \Mon(E) \to \pi_0 \cat{F}$. Assume that the functor $h_E''$ is full, and is the identity on objects. In this case, we call $E$ equipped with $h_E$ a \Def{generating graph} for $\pi_0 \cat{F}$.

Next, choose a graph morphism $s_E \colon E \to U F_0$ which is a lift of $h_E$, as in the diagram
\[
\xymatrix{
& U \cat{F}_0 \ar@{->>}[d] \\
E \ar@{-->}[ur]^-{s_E} \ar[r]_-{h_E} & U \pi_0 \cat{F}.
}
\]
Explicitly, this amounts to choosing a representative in $\cat{F}_0$ for each map $h_E(f)$ in $\pi_0 \cat{F}$. 
By adjunction, $s_E$ yields a functor $s_E' \colon \Mon(E) \to \cat{F}_0$. Since the hom-sets in the category $\cat{F}_0$ are $\F_p$-modules, we obtain by adjunction an $\F_p$-linear map
\[
s_E'' \colon \F_p \Mon(E) (A,B) \to \cat{F}(A,B)_0 
\]
for all objects $A,B$ of $\Mon(E)$, namely the vertices of $E$. Note however that $s_E''$ does \emph{not} define a functor $\F_p \Mon(E) \to \cat{F}_0$, since $\cat{F}_0$ need not be right linear. 

\begin{proposition}\label{pr:PseudoFunctor}
Let $\cat{F}$ be a left linear track category in which every morphism is $p$-torsion, equipped with a system of linearity tracks $\Ga_a^{x,y} \colon a(x+y) \Ra ax + ay$ satisfying the linearity track equations. Let $E$ be a generating graph for $\pi_0 \cat{F}$, and let $s_E'' \colon \F_p \Mon(E) \to \cat{F}_0$ be as constructed above. Let $s \colon \cat{B}_0 \to \cat{F}$ be defined as the composite
\[
\xymatrix{
\cat{B}_0 \dfn \G \Mon(E) \ar@{->>}[r] & \F_p \Mon(E) \ar[r]^-{s_E''} & \cat{F}_0 \\ 
}
\]
where the first map is the canonical quotient, induced by the quotient map $\G \surj \F_p$.
Then there exists a unique pseudo-functor $(s,\Ga) \colon \cat{B}_0 \to \cat{F}$ satisfying the following conditions.
\begin{enumerate}
\item $\Ga$ is left linear:
\[
\Ga(x+x',y) = \Ga(x,y) + \Ga(x',y).
\]
\item $\Ga(x,w) = \id^{\sq}$ if $w \in \Mon(E)$, where $\id^{\sq} \colon s(xw) \Ra s(xw)$ is the identity track.
\item  \label{eq:RightLinCorrection} $\Ga(x,y+z) = \left( \Ga(x,y) + \Ga(x,z) \right) \sq \Ga_{sx}^{sy,sz}$. 
In other words, the following diagram of tracks commutes:
\[
\xymatrix @C=4pc {
s(x) s(y+z) \ar@{=}[d] \ar@{=>}[r]^-{\Ga(x,y+z)} & s \left( x (y+z) \right) \ar@{=}[d] \\
sx (sy + sz) \ar@{=>}[d]_{\Ga_{sx}^{sy,sz}} & s(xy + xz) \ar@{=}[d] \\
(sx)(sy) + (sx)(sz) \ar@{=>}[r]_-{\Ga(x,y) + \Ga(x,z)} & s(xy) + s(xz). \\ 
}
\]
\end{enumerate}
\end{proposition}

\begin{proof}
Uniqueness follows from the fact that every morphism $y \in \cat{B}_0 = \G \Mon(E)$ is a $\G$-linear combination $y = \sum_i c_i w_i$ of words $w_i \in \Mon(E)$, in particular a finite sum of words $w_i$. 
Condition (2) determines the value of $\Ga(x,w)$ for $w \in \Mon(E)$. Applying condition (3) repeatedly then determines the value of $\Ga(x,y)$ for arbitrary $y$.

For existence, applying condition~\eqref{eq:RightLinCorrection} inductively, together with Proposition~\ref{pr:BreakSum}, yields the equality
\begin{equation} \label{eq:PseudoCorrection}
\Ga(x, \sum_{i=1}^k y_i) = \left( \sum_{i=1}^k \Ga(x,y_i) \right) \sq \Ga_{sx}^{sy_1, sy_2, \ldots, sy_k}.
\end{equation}
In other words, the following diagram of tracks commutes:
\[
\xymatrix @C=4pc {
s(x) s( \sum_{i=1}^k y_i ) \ar@{=}[d] \ar@{=>}[r]^-{\Ga(x, \sum_{i=1}^k y_i )} & s \left( x (\sum_{i=1}^k y_i) \right) \ar@{=}[d] \\
sx (\sum_{i=1}^k sy_i) \ar@{=>}[d]_{\Ga_{sx}^{sy_1, sy_2, \ldots, sy_k}} & s( \sum_{i=1}^k xy_i ) \ar@{=}[d] \\
\sum_{i=1}^k (sx)(sy_i) \ar@{=>}[r]_-{\sum_{i=1}^k \Ga(x,y_i)} & \sum_{i=1}^k s(xy_i). \\ 
}
\]
The formula~\eqref{eq:PseudoCorrection} does not depend on the ordering of the terms $y = \sum_{i=1}^k y_i$, by the symmetry equation Proposition~\ref{pr:IteratedLinTrackEq}~\eqref{item:IteratedSymm}. 
Let us check that the formula is well-defined over the ground ring $\G$. For an integer $k \in \Z$, consider the morphism in $\cat{B}_0$ given by the sum $k \cdot y = y + \ldots + y$. The diagram above specializes to
\[
\xymatrix @C=4pc {
s(x) s( k \cdot y ) \ar@{=}[d] \ar@{=>}[r]^-{\Ga(x, k \cdot y )} & s \left( x (k \cdot y) \right) \ar@{=}[d] \\
sx (k \cdot sy) \ar@{=>}[d]_{\Ga_{sx}^{sy, sy, \ldots, sy}} & s( k \cdot xy ) \ar@{=}[d] \\
k \cdot (sx)(sy) \ar@{=>}[r]_-{k \cdot \Ga(x,y)} & k \cdot s(xy). \\ 
}
\]
The corresponding equation is:
\begin{align*}
\Ga(x, k \cdot y) &= \left( k \cdot \Ga(x,y) \right) \sq \Ga_{sx}^{sy, sy, \ldots, sy} \\
&= \left( k \cdot \Ga(x,y) \right) \sq \Ga_{sx}^{1, 1, \ldots, 1} (sy) \\
&= \left( k \cdot \Ga(x,y) \right) \sq \Ga(k)_{sx} (sy).
\end{align*}
If $p^2 \vert k$ holds, then this equation of tracks yields:
\begin{align*}
\Ga(x, k \cdot y) &= \left( k \cdot \Ga(x,y) \right) \sq \Ga(k)_{sx} (sy) \\
&= \id^{\sq}_0 \sq \id^{\sq}_0 \\
&= \id^{\sq}_0
\end{align*}
where we used Lemma~\ref{lem:GammaP2}. For the left variable $x$ in $\Ga(x,y)$, a single factor of $p$ is enough:
\[
\Ga(p \cdot x, y) = p \cdot \Ga(x,y) = \id^{\sq}_0.
\]
Thus, given $\G$-linear combinations $x = \sum_i c_i x_i$ and $y = \sum_j d_j y_j$ in $\cat{B}_0 = \G \Mon(E)$, lift those to $\Z$-linear combinations $\sum_i c'_i x_i$ and $y = \sum_j d'_j y_j$ and \emph{define} $\Ga$ by the following formulas.
\begin{enumerate}
\item For arbitrary $x,y \in \G \Mon(E)$:
\[
\Ga(x,y) = \Ga \left( \sum_i c_i x_i, \sum_j d_j y_j \right) \dfn \sum_i c'_i \Ga \left( x_i, \sum_j d_j y_j \right)
\]
which does not depend on the lifts of the scalars $c_i \in \G$ to $c'_i \in \Z$.
\item \label{item:LinCombin} When $x = x_i \in \Mon(E)$ is a single word:
\[
\Ga \left( x_i, \sum_j d_j y_j \right) \dfn \left( \sum_{j=1}^k \Ga(x_i, d_j y_j) \right) \sq \Ga_{sx_i}^{d'_1 (sy_1), d'_2 (sy_2), \ldots, d'_k (s y_k)}
\]
which does not depend on the lifts of the scalars $d_j \in \G$ to $d'_j \in \Z$. 
\item \label{item:ScalarMultiple} When moreover $y$ is a scalar multiple of a word $y_j \in \Mon(E)$:
\[
\Ga(x_i, d_j y_j) \dfn \left( d'_j \cdot \Ga(x_i,y_j) \right) \sq \Ga(d'_j)_{sx_i} (sy_j)
\]
which again does not depend on the lifts $d'_j \in \Z$. The result of the previous two steps does not depend on the way to write $\sum_j d_j y_j$ as a $\G$-linear combination, by Proposition~\ref{pr:BreakSum}; for example: $\Ga(x,2y_1 + 5y_1) = \Ga(x,7y_1)$.
\item For single words $x_i, y_j \in \Mon(E)$, define:
\[
\Ga(x_i,y_j) = \id^{\sq}.
\]
\end{enumerate}
Then $\Ga(x,y)$ is well-defined, and one readily checks that $\Ga$ satisfies the three conditions in the statement. 

A straightforward (if tedious) verification shows that $(s,\Ga) \colon \G \Mon(E) \to \cat{F}$ satisfies the composition equation~\eqref{eq:CoherentIso} of a pseudo-functor \cite{Baues06}*{Theorem~5.2.3}. Also, $(s,\Ga)$ satisfies the strict unital condition. The equations $s(1_A) = 1_{sA} = 1_A$ and $\Ga(x,1) = \id^{\sq}_{sx}$ hold by construction, while the equation $\Ga(1,y) = \id^{\sq}_{sy}$ follows from the unital equation for the linearity tracks $\Ga_1^{x,y} = \id^{\sq}_{x+y}$.
\end{proof}

The proof also yields an analogous statement over $\Z$ instead of $\F_p$.

\begin{proposition} \label{pr:PseudoFunctorZ}
Let $\cat{F}$ be a left linear track category equipped with a system of linearity tracks $\Ga_a^{x,y} \colon a(x+y) \Ra ax + ay$ satisfying the linearity track equations. Let $E$ be a generating graph for $\pi_0 \cat{F}$, and let $s_E'' \colon \Z \Mon(E) \to \cat{F}_0$ be as constructed above. Denote $\cat{B}_0 \dfn \Z \Mon(E)$. 
Then there exists a unique pseudo-functor $(s,\Ga) \colon \cat{B}_0 \to \cat{F}$ satisfying the conditions listed in Proposition~\ref{pr:PseudoFunctor}.
\end{proposition}

\begin{proof}
The proof is similar to that of Proposition~\ref{pr:PseudoFunctor}, with the following changes.
Condition~\eqref{eq:RightLinCorrection} specialized to $z= -y$ yields the commutative diagram of tracks:
\[
\xymatrix @C=4.8pc {
0 = s(x) s(y + (-y)) \ar@{=}[d] \ar@{=>}[r]^-{\Ga(x,y+(-y))} & s \left( x (y + (-y)) \right) = 0 \ar@{=}[d] \\
s(x) \left( sy + s(-y) \right) \ar@{=>}[d]_{\Ga_{sx}^{sy, s(-y)}} & s \left( xy + x (-y) \right) \ar@{=}[d] \\  
(sx)(sy) + (sx)(s(-y)) \ar@{=>}[r]_-{\Ga(x,y) + \Ga(x,-y)} & s(xy) + s(x(-y)) = 0
}
\]
which in turn yields the commutative diagram
\[
\xymatrix{
**[l] (sx) (- (sy)) = s(x) s(-y) \ar@{=>}[d]_{\Ga(-1)_{sx}(sy)} \ar@{=>}[r]^-{\Ga(x,-y)} & s (x(-y)) \ar@{=}[d] \\
- (sx)(sy) \ar@{=>}[r]_-{-\Ga(x,y)} & -s(xy). \\
}
\]
In particular, the track $\Ga(x,-y)$ is determined by $\Ga(x,y)$, which proves uniqueness of $\Ga$.

In the explicit construction of $\Ga$, the result of steps~\eqref{item:LinCombin} and \eqref{item:ScalarMultiple} does not depend on the way to write $\sum_j d_j y_j$ as a $\Z$-linear combination, by Lemma~\ref{lem:SumNegatives}; for example: \linebreak $\Ga(x,-2y_1 + 5y_1) = \Ga(x,3y_1)$.
\end{proof}

\section{Pseudo-functors and associated action} \label{sec:PseudoFunctorAction}

In this section, let $\cat{T}$ be a pointed track category and let $(s,\Ga) \colon \cat{B}_0 \to \cat{T}$ be a pseudo-functor as in Definition~\ref{def:PseudoFunctor}. We will construct an action associated to a pseudo-functor, as in \cite{Baues06}*{\S 5.3}.

\subsection{The multiplicative structure}

\begin{notation} \label{no:Bullet}
Given $x,y \in \cat{B}_0$ and a track $a \colon sx \Ra 0$ in $\cat{T}_1$, define operations
\[
\begin{cases}
y \bu a = (sy \ot a) \sq \Ga(y,x)^{\inv} \in \cat{T}_1 \\
a \bu y = (a \ot sy) \sq \Ga(x,y)^{\inv} \in \cat{T}_1 \\
\end{cases}
\]
as illustrated in the commutative diagrams of tracks
\[
\xymatrix{
(sy)(sx) \ar@{=>}[d]_{\Ga(y,x)} \ar@{=>}[dr]^{sy \ot a} & \\
s(yx) \ar@{=>}[r]_{y \bu a} & 0 \\
}
\quad \quad
\xymatrix{
(sx)(sy) \ar@{=>}[d]_{\Ga(x,y)} \ar@{=>}[dr]^{a \ot sy} & \\
s(xy) \ar@{=>}[r]_{a \bu y} & 0. \\
}
\]
\end{notation}

\begin{definition} \label{def:Candidate}
Let $\cat{B}_1$ be the pullback in the diagram of (pointed) sets
\[
\xymatrix{
\cat{B}_1 \pb \ar[d]_{d} \ar[r]^{s} & \ker \de_1 \ar[d]^{\de_0} \\
\cat{B}_0 \ar[r]^{s} & \cat{T}_0 \\
}
\]
with $\ker \de_1 = \left\{ a \in \cat{T}_1 \mid \del_1 a = 0 \right\}$. Explicitly, elements of $\cat{B}_1$ are pairs $(a,x) \in \ker \de_1 \x \cat{B}_0$ satisfying $\de_0 a = sx$. Define the left and right \Def{$\ot$-action} of $\cat{B}_0$ on $\cat{B}_1$ by the formulas
\begin{equation} \label{eq:TensorAction}
\begin{cases}
y \ot (a,x) = \left( y \bu a, yx \right) \in \cat{B}_1 \\ 
(a,x) \ot y = \left( a \bu y, xy \right) \in \cat{B}_1 \\
\end{cases}
\end{equation}
for $(a,x) \in \cat{B}_1$ and $y \in \cat{B}_0$, using Notation~\ref{no:Bullet}.
\end{definition}

\begin{remark}
If one denotes the pair $(a,x) \in \cat{B}_1 = \cat{T}_1 \x_{\cat{F}_0} \cat{B}_0$ as a single symbol $\al = (a,x)$, then by definition we have $s \al = a$, $d \al = x$, and the formulas~\eqref{eq:TensorAction} can be rewritten as
\begin{equation} \label{eq:TensorActionSingle}
\begin{cases}
y \ot \al = \left( y \bu (s \al), y (d \al) \right) \in \cat{B}_1 \\ 
\al \ot y = \left( (s \al) \bu y, (d \al) y \right) \in \cat{B}_1. \\
\end{cases}
\end{equation}
\end{remark}

\begin{proposition} \label{pr:TensorAssoc}
\begin{enumerate}
\item \label{item:TensorAssoc} The $\ot$-action on $\cat{B}$ is associative. Explicitly, for $(a,x) \in \cat{B}_1$ and $y,z \in \cat{B}_0$, the following equations hold:
\[
\begin{cases}
\left( (a,x) \ot y \right) \ot z = (a,x) \ot (yz) \\
\left( y \ot (a,x) \right) \ot z = y \ot \left( (a,x) \ot z \right) \\
y \ot \left( z \ot (a,x) \right) = (yz) \ot (a,x). \\
\end{cases}
\]
\item \label{item:Unital} The $\ot$-action on $\cat{B}$ is unital, i.e., satisfies $x \ot 1 = x = 1 \ot x$ for all $x \in \cat{B}$.
\item \label{item:Leibniz} $\cat{B}$ satisfies the Leibniz rule. Explicitly, given $(a,x), (b,y) \in \cat{B}_1$, the following equation holds in $\cat{B}_1$:
\[
\left( d (a,x) \right) \ot (b,y) = (a,x) \ot \left( d(b,y) \right).
\]
Given $(a,x) \in \cat{B}_1$ and $y \in \cat{B}_0$, the following equations hold in $\cat{B}_0$:
\[
\begin{cases}
d \left( (a,x) \ot y \right) = d \left( (a,x) \right) \ot y \\
d \left( y \ot (a,x) \right) = y \ot d \left( (a,x) \right).\\
\end{cases}
\]
\end{enumerate}
\end{proposition}

\begin{proof}
\eqref{item:TensorAssoc} We will prove the equality
\[
\left( (a,x) \ot y \right) \ot z = (a,x) \ot (yz)
\]
as the other two cases are proved similarly. Expanding in terms of the product $\bu$, the equation becomes
\[
\left( (a \bu y) \bu z, xyz \right) = \left( a \bu (yz), xyz \right)
\]
or equivalently, $(a \bu y) \bu z = a \bu (yz)$. The factorization equation for the tracks $a \colon sx \Ra 0$ and $\Ga(y,z) \colon (sy)(sz) \Ra s(yz)$ in $\cat{T}_1$ yields the equality of tracks
\[
a \ot (sy)(sz) = \left( a \ot s(yz) \right) \sq \left( sx \ot \Ga(y,z) \right).
\]
Using the definition of $\bu$ and the coherence equation for the pseudo-functor $s$, the right-hand side becomes
\begin{align*}
\left( a \ot s(yz) \right) \sq \left( sx \ot \Ga(y,z) \right) &= a \bu (yz) \sq \Ga(x,yz) \sq \left( sx \ot \Ga(y,z) \right) \\
&= a \bu (yz) \sq \Ga(xy,z) \sq \left( \Ga(x,y) \ot sz \right)
\end{align*}
while the left-hand side becomes
\begin{align*}
a \ot (sy)(sz) &= \left( (a \bu y) \ot sz \right) \sq \left( \Ga(x,y) \ot sz \right) \\
&= (a \bu y) \bu z \sq \Ga(xy,z) \sq \left( \Ga(x,y) \ot sz \right)
\end{align*}
which yields the desired equality $a \bu (yz) = (a \bu y) \bu z$.

\eqref{item:Unital} For an element $x \in \cat{B}$ of degree $0$, the equations $x \ot 1 = x = 1 \ot x$ hold by definition. Now let $(a,x) \in \cat{B}$ be an element of degree $1$. The equations $1 \ot (a,x) = (a,x)$ and $(a,x) \ot 1$ are equivalent respectively to $1 \bu a = a$ and $a \bu 1 = a$. We have
\begin{align*}
1 \bu a &= (s1 \ot a) \sq \Ga(1,x)^{\inv} \\
&= a \sq \id^{\sq}_{sx} \\
&= a
\end{align*}
and likewise $a \bu 1 = a$.

\eqref{item:Leibniz} The second and third equations hold by definition of the $\ot$-action:
\begin{align*}
d \left( (a,x) \ot y \right) &= d \left( a \bu y, xy \right) \\
&= xy \\
d \left( (a,x) \right) \ot y &= x \ot y = xy.
\end{align*}
For the first equation, the two sides are:
\begin{align*}
\left( d (a,x) \right) \ot (b,y) &= x \ot (b,y) \\
&= (x \bu b, xy) \\
(a,x) \ot \left( d(b,y) \right) &= (a,x) \ot y \\
&= (a \bu y, xy)
\end{align*}
so that the equation in $\cat{B}_1$ is equivalent to the equation $x \bu b = a \bu y$ in $\cat{T}_1$. By definition of~$\bu$, we have the equalities in $\cat{T}_1$
\[
\begin{cases}
x \bu b = \left( sx \ot b \right) \sq \Ga(x,y)^{\inv} \\
a \bu y = \left( a \ot sy \right) \sq \Ga(x,y)^{\inv} \\
\end{cases}
\]
so that the equation is equivalent to $sx \ot b = a \ot sy$. This is an instance of the factorization equation \eqref{eq:FactorizToZero}:
\begin{align*}
&(\del a) \ot b = a \ot (\del b) \\
&sx \ot b = a \ot sy
\end{align*}
which holds in any pointed track category.
\end{proof}

\subsection{The left linear case} \label{sec:LeftLinearCase}

We are interested in the situation where the pointed track category $\cat{F}$ is left linear. The pullback diagram in Definition~\ref{def:Candidate} can be rewritten as
\begin{equation} \label{eq:PullbackLeftLin}
\xymatrix{
\cat{B}_1 \pb \ar[d]_{d} \ar[r]^{s} & \cat{F}_1 \ar[d]^{\del} \\
\cat{B}_0 \ar[r]^{s} & \cat{F}_0 \\
}
\end{equation}
where $\del \colon \cat{F}_1 \to \cat{F}_0$ is a morphism of abelian groups. We focus on the case where the following assumptions hold.

\begin{assumption} \label{ass:EquivHomot}
\begin{enumerate}
\item \label{eq:Preadditive} $\cat{B}_0$ is a preadditive category.
\item \label{eq:LocallyLinear} $s \colon \cat{B}_0 \to \cat{F}_0$ is locally linear. That is, for all objects $A$ and $B$ of $\cat{B}_0$, the map $s \colon \cat{B}_0(A,B) \to \cat{F}(A,B)_0$ is a morphism of abelian groups.
\item \label{eq:LocallyQuasiIso} The functor $\cat{B}_0 \ral{s} \cat{F}_0 \surj H_0 \cat{F} = \coker \del$ is full, i.e., each map $s \colon \cat{B}_0(A,B) \to H_0 \cat{F}(A,B)$ is surjective.
\item The functor $\cat{B}_0 \ral{s} \cat{F}_0 \surj H_0 \cat{F}$ is essentially surjective.
\end{enumerate}
\end{assumption}

The first two assumptions ensure that $d \colon \cat{B}_1 \to \cat{B}_0$ is a homomorphism, and Diagram~\eqref{eq:PullbackLeftLin} defines a chain map $s \colon \cat{B} \to \cat{F}$. Denote by $\si \colon H_i \cat{B} \to H_i \cat{F}$ the map induced on homology, for $i =0,1$. By the third assumption, $\si$ is an isomorphism. The fourth assumption then implies that $\si \colon H_0 \cat{B} \ral{\simeq} H_0 \cat{F}$ is an equivalence of categories. 

\begin{proposition} \label{pr:1Alg}
Let $(s,\Ga) \colon \cat{B}_0 \to \cat{F}$ be a pseudo-functor where $\cat{B}_0$ is a preadditive category and $s \colon \cat{B}_0 \to \cat{F}_0$ is locally linear. Let $(s,\Ga) \colon \cat{B} \to \cat{F}$ be defined as in Diagram~\eqref{eq:PullbackLeftLin}. Then $\cat{B}$ is a $1$-truncated DG-category.
\end{proposition}

\begin{proof}
In view of Proposition~\ref{pr:TensorAssoc} and Example~\ref{ex:1Alg}, the statement amounts to the $\ot$-composition in $\cat{B}$ being right linear. Since $\cat{B}_0$ is a preadditive category, the $\ot$-composition $x \ot y$ is bilinear in the case $\deg(x) = \deg(y) = 0$.

Let us prove the case $\deg(x) = 1$, $\deg(y) = 0$. Let $(a,x) \in \cat{B}_1$, $y,y' \in \cat{B}_0$. We want to show that $(a,x) \ot (y+y') = (a,x) \ot y + (a,x) \ot y'$ holds, which is equivalent to
\[
a \bu (y + y') = a \bu y + a \bu y'.
\] 
Consider the diagram of tracks in $\cat{F}$
\[
\xymatrix{
sx \ot (s(y+y')) = sx \ot (sy + sy') \ar@{=>}[dd]_{\Ga(x,y+y')} \ar@{=>}[dr]_-{\Ga_{sx}^{sy,sy'}} \ar@/^4.5pc/@{=>}[ddrr]^-{a \ot s(y+y')} & & \\
& (sx)(sy) + (sx)(sy') \ar@{=>}[dl]_-{\Ga(x,y) + \Ga(x,y')} \ar@{=>}[dr]_-(0.3){a \ot sy + a \ot sy'} & \\
s(x(y+y')) = s(xy) + s(xy') \ar@{=>}[rr]_-(0.35){a \bu (y+y')}^-(0.35){a \bu y + a \bu y'} & & 0 \\
}
\]
The left triangle commutes, by Proposition~\ref{pr:PseudoFunctor}~\eqref{eq:RightLinCorrection}. The top right triangle commutes, by the naturality Equation~\ref{def:LinTrackEq}~\eqref{item:NaturalA} and $\Ga_0^{sy,sy'} = \id^{\sq}_0$.  The track $a \bu (y+y')$ makes the outer triangle commute, while $a \bu y + a \bu y'$ makes the bottom triangle commute, proving the equality.

The case $\deg(x) = 0$, $\deg(y) = 1$ is similar, using the naturality Equation~\ref{def:LinTrackEq} \eqref{item:NaturalXY} and $\Ga_{a}^{0,0} = \id^{\sq}_0$. 
\end{proof}
Compare with \cite{Baues06}*{Theorem~A.15, Theorem~5.3.5}.

\begin{example}
Consider the Eilenberg--MacLane mapping theory $\cat{EM}$ and the left linear track category $\Pi_{1} \cat{EM}$. Then Proposition~\ref{pr:1Alg} yields $s \colon \cat{B} \to \Pi_{1} \cat{EM}$. This $1$-truncated DG-category $\cat{B}$ over $\G$ is called the \Def{DG-category of secondary cohomology operations}.

Given a spectrum $X$ used as distinguished object of $\Pi_{1} \cat{EM} \{X\}$ --- i.e., where we allow maps out of $X$ but never into $X$ --- Proposition~\ref{pr:1Alg} yields $\cat{B} \{X\}$ where $X$ is still a distinguished object. The $1$-truncated DG-module over $\cat{B}$
\[
H^{\strict}(X) \dfn \cat{B} \{X\} (X,-) \colon \cat{B} \to \Ch_{\leq 1}
\]
is called the \Def{strictified secondary cohomology} of $X$.
\end{example}

In \cref{sec:ChoiceGraph}, we will describe an explicit choice of generating graph $E$ which is adapted to this case.

\begin{warning}
What was called the \emph{secondary Steenrod algebra} in \cite{Baues06}*{\S 2.5} is the groupoid-enriched full subcategory of $\Pi_{1} \cat{EM}$ on the objects $\left\{ K_n \mid n \in \Z \right\}$. Likewise, what was called \emph{strictification of the secondary Steenrod algebra} in \cite{Baues06}*{Definition~5.5.2} is the $\Ch_{\leq 1}$-enriched full subcategory of our $\cat{B}$ on the objects $\left\{ K_n \mid n \in \Z \right\}$.
\end{warning}

\section{Strictification via pseudo-functors}\label{sec:Strictification}

In this section, we show how a pseudo-functor can induce a strictification, relying on a $2$-categorical observation due to Lack \cite{Lack02}, which was kindly pointed out to us by Emily Riehl. The construction we will describe is also found in \cite{Lack04}*{\S 1} and \cite{Gurski13}*{\S 4}. Let us recall some terminology.

\begin{definition}
A track functor $F \colon \cat{S} \to \cat{T}$ between track categories is called a \Def{Dwyer--Kan equivalence}, or \Def{DK-equivalence} for short, if it satisfies the following conditions.
\begin{itemize}
\item For all objects $A,B$ of $\cat{S}$, the induced map of groupoids $F \colon \cat{S}(A,B) \to \cat{T}(FA,FB)$ is an equivalence.
\item The induced functor on homotopy categories $\pi_0 F \colon \pi_0 \cat{S} \to \pi_0 \cat{T}$ is an equivalence of categories.
\end{itemize}
Track categories $\cat{S}$ and $\cat{T}$ are said to be \Def{weakly equivalent} if there is a zigzag of DK-equivalences between them.

A \Def{pseudo-DK-equivalence} $F \colon \cat{S} \to \cat{T}$ between track categories is a pseudo-functor satisfying the conditions listed above.
\end{definition}

\begin{lemma}\label{lem:StrictificationPseudo}
Let $F \colon \cat{S} \to \cat{T}$ be a pseudo-DK-equivalence between track categories. Then $F$ induces a zigzag of DK-equivalences between $\cat{S}$ and $\cat{T}$.
\end{lemma}

\begin{proof}
Let us recall the adjunction described in \cite{Lack02}*{Proposition~4.2}. 
Consider the track category $\cat{S}'$ with the same objects as $\cat{S}$ and whose morphisms are words in $\cat{S}$
\[
w = [f_{\ell}] \cdots [f_2] [f_1],
\]
that is, formal composites of composable morphisms in $\cat{S}$. Tracks $[f_{\ell}] \cdots [f_1] \Ra [g_{k}] \cdots [g_1]$ in $\cat{S}'$ are defined to be the tracks $f_{\ell} \cdots f_1 \Ra g_k \cdots g_1$ in $\cat{S}$. 
(The construction $\cat{S}'$ was called the \emph{relaxation} of $\cat{S}$ in \cite{BauesJP03}*{\S 2.4}.) 
The counit $Q \colon \cat{S}' \to \cat{S}$ sends a formal composite $[f_{\ell}] \cdots [f_1]$ to the actual composite $f_{\ell} \cdots f_1$; this $Q$ is a DK-equivalence. 
The unit $P \colon \cat{S} \to \cat{S}'$ sends a morphism $f$ to the length one word $[f]$; this $P$ is a pseudo-DK-equivalence. 

Let $G \colon \cat{S}' \to \cat{T}$ be the unique track functor satisfying $GP = F$. Then $G$ is a DK-equivalence, since $F$ is a pseudo-DK-equivalence. Hence, $Q$ and $G$ provide the desired zigzag, as illustrated in the diagram
\begin{equation}\label{eq:PseudoFunctors}
\xymatrix @C=4.0pc @R=4.0pc {
\cat{S}' \ar@/_2.0pc/[d]_(0.4)Q^(0.4){\sim} \ar[dr]^-(0.6){G}_-(0.6){\sim} & \\
\cat{S} \ar@/_0.0pc/@{~>}[u]_(0.4){P}^(0.4){\sim} \ar@{~>}[r]_-{F}^-{\sim} & \cat{T}, \\
}
\end{equation}
where the squiggly arrows denote pseudo-functors.
\end{proof}

We will need a locally linear version of that construction. 
Given a locally linear track category $\cat{S}$, consider the following construction $\tild{\cat{S}}$.
\begin{itemize}
\item The objects of $\tild{\cat{S}}$ are the same as those of $\cat{S}$. 
\item $1$-morphisms in $\tild{\cat{S}}$ are formal $\Z$-linear combinations $\sum_{i=1}^n c_i w^{(i)}$ of words $w^{(i)}$ of \linebreak $1$-morphisms in $\cat{S}$, modulo the relation generated by relations of the form \linebreak $[f + g] = [f] + [g]$ for $1$-morphisms $f,g \colon A \to B$ in $\cat{S}$.
\item $2$-morphisms in $\tild{\cat{S}}$ between formal linear combinations of words
\[
\al \colon \sum_{i=1}^m c_i [f^{(i)}_{k_i}] \cdots [f^{(i)}_{2}][f^{(i)}_{1}] \Ra \sum_{j=1}^n d_j [g^{(j)}_{\ell_j}] \cdots [g^{(j)}_{2}][g^{(j)}_{1}]
\]
are the $2$-morphisms in $\cat{S}$ between the corresponding sums of composites computed in $\cat{S}$, that is:
\[
\al \colon \sum_{i=1}^m c_i f^{(i)}_{k_i} \cdots f^{(i)}_{2} f^{(i)}_{1} \Ra \sum_{j=1}^n d_j g^{(j)}_{\ell_j} \cdots g^{(j)}_{2} g^{(j)}_{1}.
\]
\end{itemize}
In general, this construction does \emph{not} make $\tild{\cat{S}}$ into a $2$-category, since the $2$-morphisms cannot be horizontally composed. However, if the locally linear track category $\cat{S}$ is bilinear to begin with, then this construction makes $\tild{\cat{S}}$ into a track category, itself also bilinear.

\begin{lemma}\label{lem:StrictificationPseudoLin}
Let $\cat{S}$ be a bilinear track category and $\tild{\cat{S}}$ the bilinear track category described above.
\begin{enumerate}
\item Consider the assignment $\tild{Q} \colon \tild{\cat{S}} \to \cat{S}$ which is the identity on objects and sends a formal linear combination of composable words to the corresponding sum of composites. Then $\tild{Q}$ is a (strict) track functor, locally linear, and moreover a DK-equivalence.
\item Consider the assignment $\tild{P} \colon \cat{S} \to \tild{\cat{S}}$ which is the identity on objects and sends a $1$-morphism $f \colon A \to B$ to the single term with a length one word $1 [f] \colon A \to B$. Then $\tild{P}$ is a canonically a pseudo-functor, locally linear, and moreover a pseudo-DK-equivalence.
\item \label{item:Pseudo} Let $\cat{T}$ be a locally linear track category and $F \colon \cat{S} \to \cat{T}$ a locally linear pseudo-functor. Then there exists a unique locally linear (strict) track functor $G \colon \tild{\cat{S}} \to \cat{T}$ satisfying $G \circ \tild{P} = F$, cf.\ Diagram~\eqref{eq:PseudoFunctors}. 
\end{enumerate}
\end{lemma}

\begin{corollary}\label{cor:StrictificationPseudoLin}
Let $F \colon \cat{S} \to \cat{T}$ be a locally linear pseudo-DK-equivalence between locally linear track categories, where $\cat{S}$ is moreover bilinear. Then $F$ induces a zigzag of locally linear DK-equivalences between $\cat{S}$ and $\cat{T}$.
\end{corollary}

We now have all the ingredients to prove the main theorem.

\begin{theorem}\label{thm:Strictification}
Let $\cat{F}$ be a left linear track category admitting linearity tracks $\Ga_a^{x,y}$ that satisfy the linearity track equations (Definition~\ref{def:LinTrackEq}). Then $\cat{F}$ is weakly equivalent to a $1$-truncated DG-category.

If moreover every morphism in $\cat{F}$ is $p$-torsion, then $\cat{F}$ is weakly equivalent to a $1$-truncated DG-category over $\Z/p^2$.
\end{theorem}

\begin{proof}
Proposition~\ref{pr:PseudoFunctorZ} (or \ref{pr:PseudoFunctor} in the $p$-torsion case) yields a pseudo-functor \linebreak $(s,\Ga) \colon \cat{B}_0 \to \cat{F}$ which satisfies Assumption~\ref{ass:EquivHomot}. The construction in \cref{sec:LeftLinearCase} yields a pseudo-DK-equivalence $(s,\Ga) \colon \cat{B} \to \cat{F}$, which moreover is locally linear. By Proposition~\ref{pr:1Alg}, $\cat{B}$ is a $1$-truncated DG-category. Corollary~\ref{cor:StrictificationPseudoLin} then yields the desired zigzag.
\end{proof}

\begin{corollary}\label{cor:IntegralEM}
Consider the \emph{integral} Eilenberg--MacLane mapping theory $\cat{EM}_{\Z}$ consisting of finite products of integral Eilenberg--MacLane spectra
\[
A = \Si^{n_1} H\Z \x \ldots \x \Si^{n_k} H\Z
\]
and mapping spaces between them. Then the underlying track category $\Pi_1 \cat{EM}_{\Z}$ is weakly equivalent to a $1$-truncated DG-category.
\end{corollary}

\begin{proof}
By Proposition~\ref{pr:LinTrackEq}, $\Pi_1 \cat{EM}_{\Z}$ is left linear and admits canonical linearity tracks $\Ga_a^{x,y}$. The result then follows from \cref{thm:Strictification}.
\end{proof}

\appendix

\section{The Steenrod algebra and a choice of generating graph} \label{sec:ChoiceGraph}

Consider the example described in the introduction $\cat{EM}$ and the left linear track category $\cat{F} = \Pi_{1} \cat{EM}$, in which every morphism is $p$-torsion. Recall that the objects of $\cat{F}$ are the finite products $K = \prod_i K_{n_i}$, with $K_{n} = \sh^n K_0 \simeq \Si^n H\F_p$ some convenient model for Eilenberg--MacLane spectra \cite{BauesFrankland17}*{Corollary~A.8}. 
We now describe how to produce a generating graph $E$ of $\pi_0 \cat{EM}$ and a lift $s_E \colon E \to U \cat{F}_0 = U \cat{EM}$ as in Section~\ref{sec:PseudoFunctor}. To begin, make the following choices.

\begin{enumerate}
\item Choose generators $E_{\steen} \subseteq \steen$ of the Steenrod algebra as an $\F_p$-algebra. Each $a \in \steen^n$ of degree $n$ corresponds to a homotopy class $a \colon H\F_p \to \Si^n H\F_p$.
\item For each generator $a \in E_{\steen}$, of degree $n$, choose a representing map $\tild{a} \colon K_0 \to K_n$ in $\cat{EM}$.
\end{enumerate}

The generating set $E_{\steen}$ of the Steenrod algebra yields a generating graph $E$ of $\pi_0 \cat{EM}$. Explicitly, $E \subseteq U \pi_0 \cat{EM}$ is the subgraph with the same vertices, and whose edges consist of matrices of elements in $E_{\steen}$, namely the homotopy classes
\[
f = \begin{bmatrix} a_{i,j} \end{bmatrix} \colon \prod_j \Si^{m_j} H\F_p \to \prod_i \Si^{n_i} H\F_p
\]
where each $a_{i,j} \colon \Si^{m_j} H\F_p \to \Si^{n_i} H\F_p$ is $\Si^{m_j} a_{i,j}'$ for some generator $a_{i,j}' \in E_{\steen}$. The shift
\[
\sh^{m_j} \tild{a}'_{i,j} \colon K_{m_j} \to \sh^{m_j} K_{n_i - m_j} = K_{n_i}
\]
is a map in $\cat{EM}$ representing $a_{i,j}$. For fixed $i$, consider the map in $\cat{EM}$
\[
\tild{f}_i \dfn \sum_j \sh^{m_j} \tild{a}'_{i,j} \circ \proj_j \colon \prod_j K_{m_j} \to K_{n_i}
\]
and let $\tild{f} \colon \prod_j K_{m_j} \to \prod_i K_{n_i}$ be the map in $\cat{EM}$ whose $i^{\text{th}}$ coordinate map is $\tild{f}_i$. By construction, $\tild{f}$ is a representative of the homotopy class $f$ in $\pi_0 \cat{EM}$.

Define the graph morphism $s_E \colon E \to U \cat{EM}$ as the identity on vertices and $s_E(f) = \tild{f}$ on edges. Then $s_E$ lifts the inclusion $h_E \colon E \to U \pi_0 \cat{EM}$.

There is an analogous construction given a spectrum $X$. Choose generators $E_X \subseteq H^* X$ of the cohomology of $X$ as an $\steen$-module, and a representing map $\tild{x} \colon X \to K_n$ in $\cat{EM} \{X\}$ for each generator $x \in E_X$, with $x \in H^n X$. Repeating the construction above, we obtain a generating graph $h_E \colon E \to U \pi_0 \cat{EM} \{X\}$ and a lift $s_E \colon E \to U \cat{EM} \{X\}$.

\section{Toda brackets via a strictification} \label{sec:Toda}

In this section, we explain how a strictification $\cat{B} \ral{\sim} \cat{F}$ as in \cref{thm:Strictification} can be used to compute Toda brackets in the homotopy category $H_0 \cat{F}$.

\subsection{Toda brackets in track categories}

\begin{definition} \label{def:TodaBracket}
Let $\cat{T}$ be a pointed track category and let 
\[
\xymatrix{
Y_0 & Y_1 \ar[l]_{y_1} & Y_2 \ar[l]_{y_2} & Y_3 \ar[l]_{y_3} \\
}
\]
be a diagram in $\pi_0 \cat{T}$ satisfying $y_1 y_2 = 0$, $y_2 y_3 = 0$. Choose maps $x_i$ in $\cat{T}_0$ representing $y_i$. Then there exist tracks $a$, $b$ as in the diagram
\[
\xymatrix{
Y_0 & Y_1 \ar[l]^{x_1} & Y_2 \ar[l]^{x_2} \lllowertwocell<-10>_0{a} & Y_3. \ar[l]^{x_3} \lllowertwocell<10>^0{^b} \\
}
\]
The $3$-fold \Def{Toda bracket} is the subset $\lan y_1, y_2, y_3 \ran \subseteq \pi_1 \cat{T}(Y_3,Y_0)$ of all elements in $\Aut(0) = \pi_1 \cat{T}(Y_3,Y_0)$ of the form
\begin{equation} \label{eq:Toda}
(a x_3) \sq (x_1 b)^{\inv}
\end{equation}
as above. Each such element in $\cat{T}_1$ is called a \Def{representative} of the Toda bracket $\lan y_1, y_2, y_3 \ran$.
\end{definition}

\begin{definition}
A pseudo-functor $(s,\Ga) \colon \cat{S} \to \cat{T}$ between pointed track categories is called \Def{pointed} if it satisfies $s(0) = 0$ and $\Ga(x,0) = \Ga(0,y) = \id^{\sq}_0$.
\end{definition}

For a morphism $x \in \cat{S}_0$ and a track $b \colon y \Ra y'$ in $\cat{S}_1$, the composition equation of the pseudo-functor $(s,\Ga) \colon \cat{S} \to \cat{T}$ reads
\[
s(xb) = \Ga(x,y') \sq (sx)(sb) \sq \Ga(x,y)^{\inv}
\]
and similarly for $s(ay) \in \cat{T}_1$. For tracks to zero $a \colon x \Ra 0$ and $b \colon y \Ra 0$ and a pointed pseudo-functor $(s,\Ga)$, this specializes to
\begin{equation} \label{eq:PointedPseudo}
\begin{cases}
s(xb) = (sx)(sb) \sq \Ga(x,y)^{\inv} \\
s(ay) = (sa)(sy) \sq \Ga(x,y)^{\inv} \\
\end{cases}
\end{equation}
as illustrated in the diagrams of tracks:
\[
\xymatrix{
s(xy) \ar@{=>}[r]^-{s(xb)} & 0 \\
(sx)(sy) \ar@{=>}[u]^{\Ga(x,y)} \ar@{=>}[ur]_-{(sx)(sb)} & \\
}
\quad \quad
\xymatrix{
s(xy)\ar@{=>}[r]^-{s(ay)} & 0 \\
(sx)(sy). \ar@{=>}[u]^{\Ga(x,y)} \ar@{=>}[ur]_-{(sa)(sy)} & \\
}
\]

\begin{proposition}\label{pr:PseudoToda}
Let $(s,\Ga) \colon \cat{S} \to \cat{T}$ be a pointed pseudo-functor between pointed track categories, and let $\si \colon \pi_i \cat{S}(A,B) \to \pi_i \cat{T}(sA,sB)$ denote the induced map on homotopy groups, for $i = 0,1$. Then:

\begin{enumerate}
\item $s \colon \cat{S} \to \cat{T}$ sends Toda bracket representatives in $\cat{S}$ to Toda bracket representatives in $\cat{T}$, in the following sense.

Given $y_1, y_2, y_3 \in \pi_0 \cat{S}$ represented by $x_1, x_2, x_3 \in \cat{S}_0$, with tracks $a, b \in \cat{S}_1$ of the form $a \colon x_1 x_2 \Ra 0$ and $b \colon x_2 x_3 \Ra 0$. 
Then the Toda bracket representative
\[
a x_3 \sq (x_1 b)^{\inv} \in \cat{S}_1 
\]
of $\lan y_1, y_2, y_3 \ran \subseteq \pi_1 \cat{S}$ is sent by $s \colon \cat{S}_1 \to \cat{T}_1$ to a Toda bracket representative
\[
s \left( a x_3 \sq (x_1 b)^{\inv} \right) = a' (sx_3) \sq \left( (sx_1) b' \right)^{\inv} \in \cat{T}_1
\]
of $\lan \si y_1, \si y_2, \si y_3 \ran \subseteq \pi_1 \cat{T}$, for some $a', b' \in \cat{T}_1$.

In particular, the following inclusion holds in $\pi_1 \cat{T}$:
\[
\si \lan y_1, y_2, y_3 \ran \subseteq \lan \si y_1, \si y_2, \si y_3 \ran.
\]

\item If moreover $s \colon \cat{S} \to \cat{T}$ induces isomorphisms $\si \colon \pi_i \cat{S}(A,B) \ral{\cong} \pi_i \cat{T}(sA,sB)$ for \linebreak $i=0,1$ and all objects $A$ and $B$ of $\cat{S}$, then the following subsets of $\pi_1 \cat{T}$ are equal: 
\[
\si \lan y_1, y_2, y_3 \ran = \lan \si y_1, \si y_2, \si y_3 \ran.
\]
\end{enumerate}
\end{proposition}

\begin{proof}
(1) Take the track $a' \dfn sa \sq \Ga(x_1,x_2)$ in $\cat{T}_1$, as illustrated here:
\[
\xymatrix{
(sx_1)(sx_2) \ar@{=>}[r]^-{\Ga(x_1,x_2)} & s(x_1 x_2) \ar@{=>}[r]^-{sa} & 0 \\ 
}
\]
and likewise $b' \dfn sb \sq \Ga(x_2,x_3)$. 
We claim the equality in $\cat{T}_1$
\begin{equation} \label{eq:PseudoToda}
s \left( a x_3 \sq (x_1 b)^{\inv} \right) = a' (sx_3) \sq \left( (sx_1) b' \right)^{\inv}.
\end{equation}
Starting from the left-hand side, we find:
\begin{flalign*}
&& &s \left( a x_3 \sq (x_1 b)^{\inv} \right) && \\
&& =&s(a x_3) \sq s(x_1 b)^{\inv} && \\
&& =&(sa)(sx_3) \sq \Ga(x_1 x_2, x_3)^{\inv} \sq \Ga(x_1, x_2 x_3) \sq \left( (sx_1)(sb) \right)^{\inv} &&\text{by Equation } \eqref{eq:PointedPseudo} \\
&& =&(sa)(sx_3) \sq \Ga(x_1,x_2)(sx_3) \sq \left( (sx_1) \Ga(x_2,x_3) \right)^{\inv} \sq \left( (sx_1)(sb) \right)^{\inv} &&\text{by the associativity equation} \\
&& =&a' (sx_3) \sq \left( (sx_1) b' \right)^{\inv} &&\text{by definition of } a' \\
\end{flalign*}
as claimed. 
The equations are illustrated in the commutative diagram of tracks in $\cat{T}_1$:
\[
\xymatrix{
0 & & (sx_1)(sx_2)(sx_3) \ar@{=>}[ll]_-{a' (sx_3)} \ar@{=>}[dl]^-*!/^0.8ex/{\labelstyle \Ga(x_1, x_2) (sx_3)} \ar@{=>}[dr]_-{(sx_1)\Ga(x_2, x_3)} \ar@{=>}[rr]^-{(sx_1) b'} & & 0. \\
& s(x_1 x_2) (sx_3) \ar@{=>}[ul]_-{(sa)(sx_3)} \ar@{=>}[dr]^-{\Ga(x_1 x_2, x_3)} & & (sx_1) s(x_2 x_3) \ar@{=>}[dl]_-{\Ga(x_1, x_2 x_3)} \ar@{=>}[ur]^-{(sx_1)(sb)} & \\
& & s(x_1 x_2 x_3) \ar@{=>}@/^4pc/[uull]^-{s(a x_3)} \ar@{=>}@/_4pc/[uurr]_-{s(x_1 b)} & & \\
}
\]

(2) Represent the maps $\si y_i \in \pi_0 \cat{T}$ by $sx_i \in \cat{T}_0$. Consider a Toda bracket representative
\[
\te = a' (sx_3) \sq \left( (sx_1) b' \right)^{\inv} \in \lan \si y_1, \si y_2, \si y_3 \ran
\]
The existence of the track $a' \colon (sx_1)(sx_2) \Ra 0$ ensures that $x_1 x_2$ is also nullhomotopic, since the map of groupoids $s \colon \cat{S}(Y_2,Y_0) \to \cat{T}(sY_2, sY_0)$ induces a bijection on $\pi_0$. 
Likewise, $b' \colon (sx_2)(sx_3) \Ra 0$ ensures $x_2 x_3 \simeq 0$. 
Moreover, $s \colon \cat{S}(Y_2,Y_0) \to \cat{T}(sY_2, sY_0)$ induces bijections on sets of tracks to $0$, which are torsors for $\pi_1 = \Aut(0)$. Hence, there exist tracks $a \colon x_1 x_2 \Ra 0$ and $b \colon x_2 x_3 \Ra 0$ in $\cat{S}_1$ satisfying
\[
\begin{cases}
sa = a' \sq \Ga(x_1,x_2)^{\inv} \\
sb = b' \sq \Ga(x_2,x_3)^{\inv}. \\
\end{cases}
\]
By Equation~\eqref{eq:PseudoToda}, the track
\[
\te = s \left( a x_3 \sq (x_1 b)^{\inv} \right)
\]
lies in the image of the restriction $\si \colon \lan y_1, y_2, y_3 \ran \to \pi_1 \cat{T}$.
\end{proof}

\subsection{Massey products in DG-categories}

\begin{definition}\label{def:MasseyProd}
Let $\cat{B}$ be a $1$-truncated DG-category and let 
\[
\xymatrix{
Y_0 & Y_1 \ar[l]_{y_1} & Y_2 \ar[l]_{y_2} & Y_3 \ar[l]_{y_3} \\
}
\]
be a diagram in $H_0 \cat{B}$ satisfying $y_1 y_2 = 0$, $y_2 y_3 = 0$. Choose maps $x_i$ in $\cat{B}_0$ representing $y_i$. Since $x_1 x_2$ and $x_2 x_3$ are zero in homology, there exist elements $a, b \in \cat{B}_1$ satisfying $\del a = x_1 x_2$ and $\del b = x_2 x_3$. The element $a x_3 - x_1 b \in \cat{B}(Y_3,Y_0)_1$ is a cycle:
\begin{align*}
\del (a x_3 - x_1 b) &= \del (a x_3) - \del(x_1 b) \\
&= (\del a) x_3 - x_1 (\del b) \\
&= (x_1 x_2) x_3 - x_1 (x_2 x_3) \\
&= 0.
\end{align*}
The $3$-fold \Def{Massey product} is the subset $\lan y_1, y_2, y_3 \ran \subseteq H_1 \cat{B}(Y_3,Y_0)$ of all elements in $\ker \del = H_1 \cat{B}(Y_3,Y_0)$ of the form
\[
a x_3 - x_1 b
\]
as above. Each such cycle in $\cat{B}_1$ is called a \Def{representative} of the Massey product $\lan y_1, y_2, y_3 \ran$.
\end{definition}

\begin{remark}\label{rem:MasseyProdGeneral}
Definition~\ref{def:MasseyProd} works more generally in a locally linear track category $\cat{F}$, where the mapping groupoids are viewed as $1$-truncated chain complexes $\cat{F}(A,B)$ via the equivalence from Proposition~\ref{pr:AbGpObjGpd}. One could instead work in the underlying pointed track category $\denorm \cat{F}$. Via the correspondence $\pi_i \denorm \left( \cat{F}(A,B) \right) \cong H_i \cat{F}(A,B)$, the Toda bracket in Definition~\ref{def:TodaBracket} corresponds to the Massey product, so that there is no ambiguity in the notation $\lan y_1, y_2, y_3 \ran$.
\end{remark}

Specializing Proposition~\ref{pr:PseudoToda} to the setup of \cref{sec:LeftLinearCase} yields the following.

\begin{corollary}\label{cor:MasseyToToda}
Let $\cat{F}$ be a left linear track category. Let $(s,\Ga) \colon \cat{B}_0 \to \cat{F}$ be a pseudo-functor where $\cat{B}_0$ is a preadditive category and $s \colon \cat{B}_0 \to \cat{F}_0$ is locally linear. Let $s \colon \cat{B} \to \cat{F}$ be the pseudo-functor as in Proposition~\ref{pr:1Alg}. 
Then:
\begin{enumerate}
\item $s \colon \cat{B} \to \cat{F}$ sends Massey product representatives in $\cat{B}$ to Toda bracket representatives in $\cat{F}$, in the following sense. 

Given $y_1, y_2, y_3 \in H_0 \cat{B}$ represented by $x_1, x_2, x_3 \in \cat{B}_0$, with $a, b \in \cat{B}_1$ satisfying $d a = x_1 x_2$ and $d b = x_2 x_3$. 
Then the Massey product representative
\[
a x_3 - x_1 b \in \cat{B}_1
\]
of $\lan y_1, y_2, y_3 \ran \subseteq H_1 \cat{B}$ is sent by $s \colon \cat{B}_1 \to \cat{F}_1$ to a Toda bracket representative
\[
s \left( a x_3 - x_1 b \right) = a' (sx_3) - (sx_1) b' \in \cat{F}_1
\]
of $\lan \si y_1, \si y_2, \si y_3 \ran \subseteq H_1 \cat{F}$, for some $a', b' \in \cat{F}_1$.

In particular, the following inclusion holds in $H_1 \cat{F}$:
\[
\si \lan y_1, y_2, y_3 \ran \subseteq \lan \si y_1, \si y_2, \si y_3 \ran.
\]
\item If moreover the functor $\cat{B}_0 \ral{s} \cat{F}_0 \surj H_0 \cat{F}$ is full (so that $s \colon \cat{B} \to \cat{F}$ is locally a quasi-isomorphism), then the following subsets of $H_1 \cat{F}$ are equal: 
\[
\si \lan y_1, y_2, y_3 \ran = \lan \si y_1, \si y_2, \si y_3 \ran.
\]
\end{enumerate}
\end{corollary}

Compare with~\cite{Baues06}*{Equation~(A18), Definition~5.5.7}.

\vspace{8pt}

\end{document}